\documentclass[twoside]{article}
\usepackage{color}
\usepackage{latexsym}
\usepackage{amsmath}
\usepackage{amsfonts}
\usepackage{amssymb}
\usepackage{graphicx}
\usepackage{indentfirst}
\usepackage{makeidx}
\usepackage{url}
\usepackage{wrapfig}
\usepackage[colorlinks=true]{hyperref}
\hypersetup{urlcolor=blue, citecolor=red}

\usepackage{times,url}
\usepackage[normal]{subfigure}
\begin{document}

\newcommand{\cC}{{\cal C}}
\newcommand{\rL}{{\rm L}}
\newcommand{\rH}{{\rm H}}
\newcommand{\rN}{{\rm N}}
\newcommand{\cN}{{\cal N}}
\newcommand{\rW}{{\rm W}}

\newcommand{\disp}{\displaystyle}
\newcommand{\proof}{\par
\noindent {\sc Proof}.\quad}
\newcommand{\fin}{\hfill\mbox{$\quad{}_{\Box}$}}
\newcommand{\fineq}{\vspace{-.75cm$\fin$}\par\bigskip}
\newcommand{\fineqnum}{\vspace{-.4cm$\fin$}\par\bigskip}

\newtheorem{coro}{\bf \sffamily Corollary}
\newtheorem{theo}{\bf \sffamily Theorem}
\newtheorem{prop}{\bf \sffamily Proposition}
\newtheorem{rem}{\bf \sffamily Remark}
\newtheorem{lemma}{\bf \sffamily Lemma}
\newtheorem{exam}{\bf \sffamily Example}
\newtheorem{defi} {\bf \sffamily Definition}
\newtheorem{nota}{\bf \sffamily Notation}

\oddsidemargin=-3mm \topmargin=-7mm

\title{\Large \bfseries \sffamily Flat solutions of some non-Lipschitz autonomous semilinear equations may be stable for $\rN\geq 3$}
\author{\bfseries\sffamily Jes\'{u}s Ildefonso D\'{\i}az,  Jes\'{u}s Hern\'{a}ndez  and Yavdat Il'yasov \thanks{J.I. D\'{\i}az and J. Hern\'{a}ndez are
partially supported by the projects ref. MTM2011-26119 and MTM2014-57113 of
the DGISPI (Spain). The research of J.I.~D\'{\i}az was partially supported
by the UCM Research Group MOMAT (Ref. 910480). The research of Y.~Il'yasov was partially supported
by RFBR-14-01-00736-a
\hfil\break\indent {\sc 2010 Mathematics Subject Classification}: 35J60, 35J96, 35R35, 53C45.
}}
\date{}
\maketitle
\begin{center}
{\em To a master, Ha\"{\i}m Brezis, with admiration.}
\end{center}

\begin{abstract}
We prove that flat ground state solutions ($i.e.$ minimizing the energy and with gradient vanishing on the boundary of the domain) of the Dirichlet problem associated to some semilinear autonomous elliptic equations with a strong absorption term given by a non-Lipschitz
function are unstable for dimensions $\rN=1,2$ and they can be stable for $\rN\geq 3$ for suitable values of the involved exponents.

\end{abstract}

\section{Introduction and main results}

Let $\rN\geq 1$, and let $\Omega $ be a bounded domain in $\mathbb{R}^{\rN}$ whose boundary $\partial \Omega $ is a $\cC^{1}$-manifold. We consider the following semi-linear parabolic problem

\begin{equation}
PP(\alpha ,\beta ,\lambda ,v_{0}) \quad\left\{ 
\begin{array}{ll}
v_{t}-\Delta v+|v|^{\alpha -1}v=\lambda |v|^{\beta -1}v & \text{in }%
(0,+\infty )\times \Omega  \\ 
v=0 & \text{on }(0,+\infty )\times \partial \Omega  \\ 
v(0,x)=v_{0}(x) & \text{on }\Omega .
\end{array}%
\right.   \label{p1}
\end{equation}%
Here $\lambda $ is a positive parameter and $0<\alpha <\beta \leq 1$. Our main goal is to give some stability criteria on solutions of the associated
stationary problem%
\begin{equation}
SP(\alpha ,\beta ,\lambda )\quad \left\{ 
\begin{array}{ll}
-\Delta u+|u|^{\alpha -1}u=\lambda |u|^{\beta -1}u & \text{in }\Omega \text{,%
} \\ 
u=0 & \text{on }\partial \Omega .%
\end{array}%
\right.   \label{1}
\end{equation}%
Notice that since the diffusion-reaction balance involves the non-linear
reaction term 
$$
f(\lambda ,u):=\lambda |u|^{\beta -1}u-|u|^{\alpha -1}u
$$ 
and it is a non-Lipschitz function at zero (since $\alpha <1$ and $\beta\le 1$) important
peculiar behavior of solutions of both problems arise. For instance, that
may lead to the violation of the Hopf maximum principle on the boundary and
the existence of compactly supported solutions as well as the so called \ 
\textit{flat solutions} (sometimes also called \textit{free boundary
solutions}) which correspond to weak solutions $u$ such that 
\begin{equation}
\frac{\partial u}{\partial \nu }=0~~\mbox{on}~~\partial \Omega ,  \label{N}
\end{equation}%
where $\nu $ denotes the unit outward normal to $\partial \Omega $. \
Solutions of this kind for stationary equations with non-Lipschitz
nonlinearity have been investigated in a number of papers. The pioneering
paper in which it was proved that the solution gives rise to a free boundary
defined as the boundary of its support was due to Ha\"{\i}m Brezis \cite%
{Brezis-Uspeki} concerning multivalued non-autonomous semilinear \
equations. The semilinear case with non-Lipschitz perturbations was
considered later in \cite{Ben-Bre-Crandall} (see also \cite{Benssousan-B-F}, 
\cite{Bre-Lieb} and \cite{brNi1996}). For the case of semilinear autonomous
elliptic equations see e.g. \cite{diaz}, \cite{diaz-H-Il}, \cite%
{Diaz-Hernan-Man}, \cite{CortElgFelmer-1}, \cite{CortElgFelmer-2}, \cite%
{IlEg}, \cite{Kaper1}, \cite{Kaper2}, \cite{Serrin-Zou}, to mention only a
few. For problem (\ref{1}), the existence of radial flat solutions was first
proved by Kaper and Kwong \cite{Kaper1}. In this paper, applying shooting
methods they showed that there exists $R_{0}>0$ such that (\ref{1})
considered in the ball $B_{R_{0}}=\{x\in \mathbb{R}^{N}:~|x|\leq
R_{0}\}=\Omega $ has a radial compactly supported positive solution.
Furthermore, by the moving-plane method it was proved in \cite{Kaper2} that any classical solution $u\in \cC^{2}(\Omega )$ of (\ref{1}) is
necessarily radially symmetric if $\Omega $ is a ball. Observe that from
this it follows that the Dirichlet boundary value problem (\ref{1}) has a
compactly supported solution if $B_{R_{0}}\subseteq \Omega $.

\bigskip

In this work we study the stability of solutions of $\ $the
stationary problem $SP(\alpha ,\beta ,\lambda )$. We point out that a direct
analysis of the stability of the stationary solutions $u_{\infty }\in
\lbrack 0,+\infty )$ of the associated ODE

\begin{equation}
ODE(\alpha ,\beta ,\lambda ,v_{0})\quad \left\{ 
\begin{array}{ll}
v_{t}+|v|^{\alpha -1}v=\lambda |v|^{\beta -1}v & \text{in }(0,+\infty ) \\ 
v(0)=v_{0}, & 
\end{array}%
\right.
\end{equation}%
shows that the trivial solution $u_{\infty }\equiv 0$ is asymptotically
stable and that the nontrivial stationary solution $u_{\infty }:=\lambda
^{-1/(\beta -\alpha )}$ is unstable (see Figure 1).

\begin{figure}[th]
\center{\includegraphics[width=0.7\linewidth]{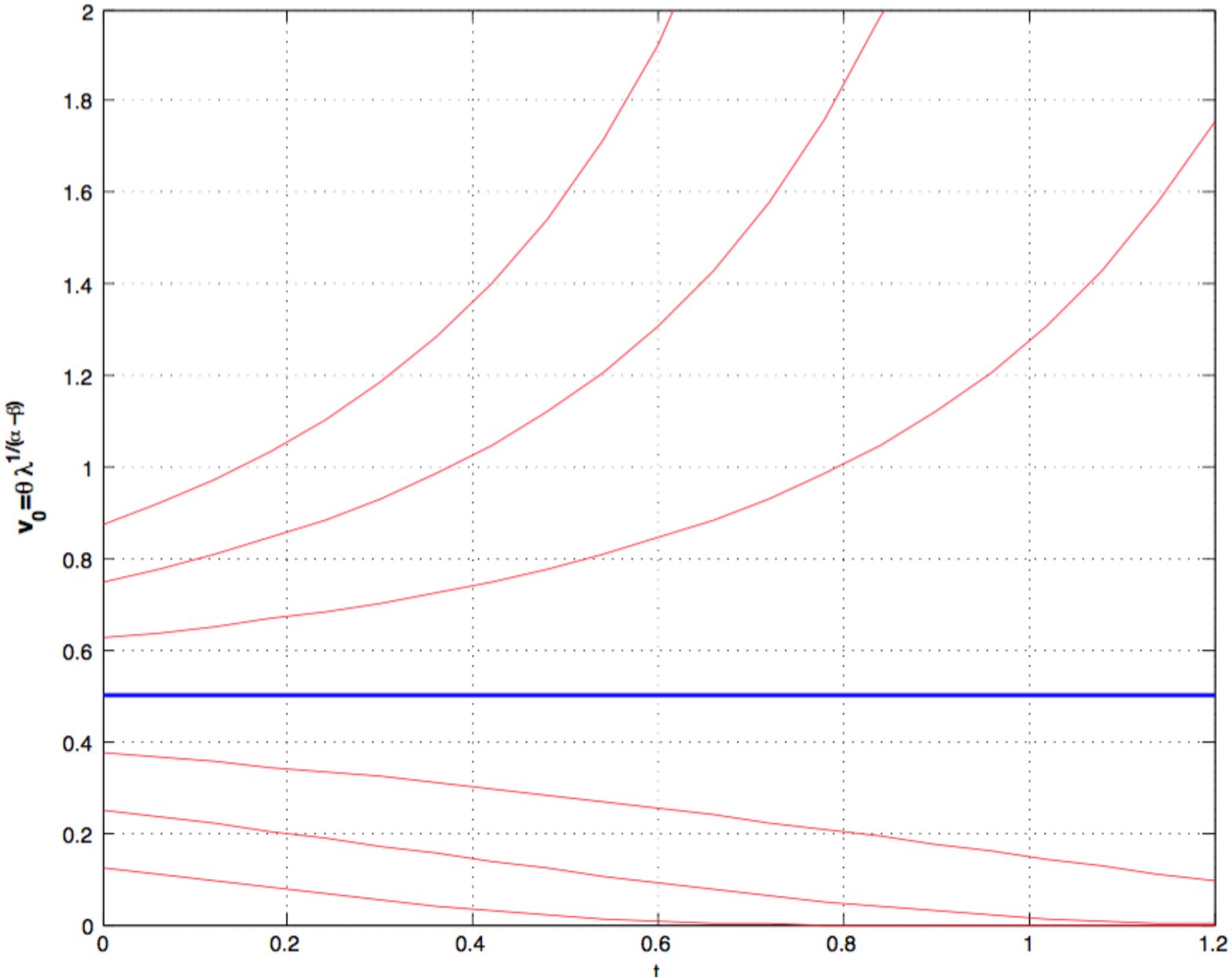}}
\caption{Paths for ${\rm ODE}(1/2,3/2,\lambda ,v_{0}).$
}
\end{figure}

Obviously the same criteria hold for the case of the semilinear problem with
Neumann boundary conditions. Nevertheless, unexpectedly, the situation is
not similar for the case of Dirichlet boundary conditions, and so, \ as the
main result of this paper will show, for dimensions $\rN\geq 3$ the nontrivial
flat solution of $SP(\alpha ,\beta ,\lambda )$ becomes stable in a certain
range of the exponents $\alpha <\beta <1$. To be more precise, our stability
study will concern \textit{ground state} solutions (also called simply \textit{ground state}) of $SP(\alpha ,\beta
,\lambda ).$ By it we mean a nonzero weak solution $u_{\lambda }$ of $%
SP(\alpha ,\beta ,\lambda )$ which satisfies 
\begin{equation*}
E_{\lambda }(u_{\lambda })\leq E_{\lambda }(w_{\lambda })
\end{equation*}%
for any nonzero weak solution $w_{\lambda }$ of $SP(\alpha ,\beta ,\lambda )$%
. Here $E_{\lambda }(u)$ is the energy functional corresponding to $%
SP(\alpha ,\beta ,\lambda )$ which is defined on the Sobolev space $%
\rH_{0}^{1}(\Omega )$ as follows 
\begin{equation*}
E_{\lambda }(u)=\frac{1}{2}\int_{\Omega }|\nabla u|^{2}\,dx+\frac{1}{\alpha
+1}\int_{\Omega }|u|^{\alpha +1}\,dx-\lambda \frac{1}{\beta +1}\int_{\Omega
}|u|^{\beta +1}\,dx.
\end{equation*}

For simplicity, we shall assume the initial value such that $v_{0}\in
\rL^{\infty }(\Omega ),~ v_{0}\ge 0$. As we shall show in Section~2, then there exists a
weak solution $v\in \cC([0,+\infty ),\rL^{2}(\Omega ))$ of $PP(\alpha ,\beta
,\lambda ,v_{0})$ satisfying $\lambda |v|^{\beta -1}v-|v|^{\alpha -1}v\in
\rL^{\infty }((0,+\infty )\times \Omega )$ \ and 
\begin{equation}
v(t)=T(t)v_{0}+\int_{0}^{t}T(t-s)(\lambda |v|^{\beta -1}v-|v|^{\alpha
-1}v)ds,  \label{wequat}
\end{equation}%
with $(T(t))_{t\geq 0}$ the heat semigroup with homogeneous Dirichlet
boundary conditions, i.e. $T(t)=e^{t(-\Delta )}$. Among some additional
regularity properties of $v$ we mention that 
\begin{equation}
v-T(t)v_{0}\in \rL^{p}(\tau ,T;\rW^{2,p}(\Omega )\cap \rW_{0}^{1,p}(\Omega ))\cap
\rW^{1,p}(\tau ,T;\rL^{p}(\Omega )),
\end{equation}%
for every $p\in (1,\infty ),$ and for any $0<\tau <T$ \big (in fact $\tau =0$ if we also
assume that $v_{0}\in \rW_{0}^{1,p}(\Omega )$\big ). In particular, $v$ satisfies
the equation $PP(\alpha ,\beta ,\lambda ,v_{0})$ for a.e. $t\in
(0,+\infty )$. Moreover, if $v(0)\in \rH_{0}^{1}(\Omega )$ then, for any $t>0$ 
\begin{equation}
\int_{0}^{t}||v_{t}(s)||_{L^{2}}^{2}ds+E_{\lambda }(v(t))\leq E_{\lambda
}(v(0)).  \label{ener1}
\end{equation}

We shall show in Section 2 that there is uniqueness of solutions
of $PP(\alpha ,\beta ,\lambda ,v_{0})$ in the class of solutions $v$ such
that 
\begin{equation}
v(t,x)\geq Cd(x)^{2/(1-\alpha )}\quad \text{in }\Omega ,\text{ for }t>0
\label{Non degenerate}
\end{equation}%
for some constant $C>0$, where $d(x):=\mathrm{dist}(x,\partial \Omega )
$ (which we shall also denote simply as $\delta _{\Omega }$). Sufficient
conditions implying this non-degeneracy property (\ref{Non degenerate}) will
be given. We also prove that if $\lambda \in \lbrack 0,\lambda _{1})$ then
the finite extinction time property is satisfied for solutions of $PP(\alpha
,\beta ,\lambda ,v_{0})$ (as in the pioneering paper \cite{Brezis-Friedman}
on multivalued semilinear parabolic problems; see also the survey \cite{Diaz
gaeta}). Moreover we shall show in Section 2 that there is a certain
resemblance between the set of solutions of $PP(\alpha ,\beta ,\lambda
,v_{0})$ and the corresponding one of the ODE problem $ODE(\alpha ,\beta
,\lambda ,v_{0})$ since: a) for any $\lambda >0$ the trivial solution $%
u\equiv 0$ of the stationary problem $SP(\alpha ,\beta ,\lambda )$ is
asymptotically stable in the sense that it attracts solutions of $PP(\alpha
,\beta ,\lambda ,v_{0})$\ for small initial data $v_{0}$ (Proposition 2.1),
and b) if $v_{0}$ is "large enough" the trajectory of the solution of $%
PP(\alpha ,\beta ,\lambda ,v_{0})$ is not non-uniformly bounded when $%
t\nearrow +\infty $ (Proposition 2.4). 

\bigskip 

Concerning the stationary problem $SP(\alpha ,\beta ,\lambda )$ we recall
that if $u\in \rH_{0}^{1}(\Omega )\cap \rL^{\infty }(\Omega )$ is a weak
stationary solution of $SP(\alpha ,\beta ,\lambda )$ then, by standard
regularity results, $u\in \rW^{2,p}(\Omega )$ for any $p\in (1,\infty )$ and
then $u\in \cC^{1,\gamma }(\overline{\Omega })$ for any $\gamma $.

In our stability study we shall use some fibrering techniques. For given $%
u\in \rH_{0}^{1}(\Omega )$, the \textit{fibrering mappings} are defined by $%
\Phi _{u}(r)=E_{\lambda }(ru)$ so that from the variational formulation of $%
SP(\alpha ,\beta ,\lambda )$ we know that $\Phi _{u}^{\prime }(r)=0$ where
we use the notation 
\begin{equation*}
\Phi _{u}^{\prime }(r)=\frac{\partial }{\partial r}E_{\lambda }(ru).
\end{equation*}%
\begin{figure}[th]
\center{\includegraphics[width=0.55\linewidth]{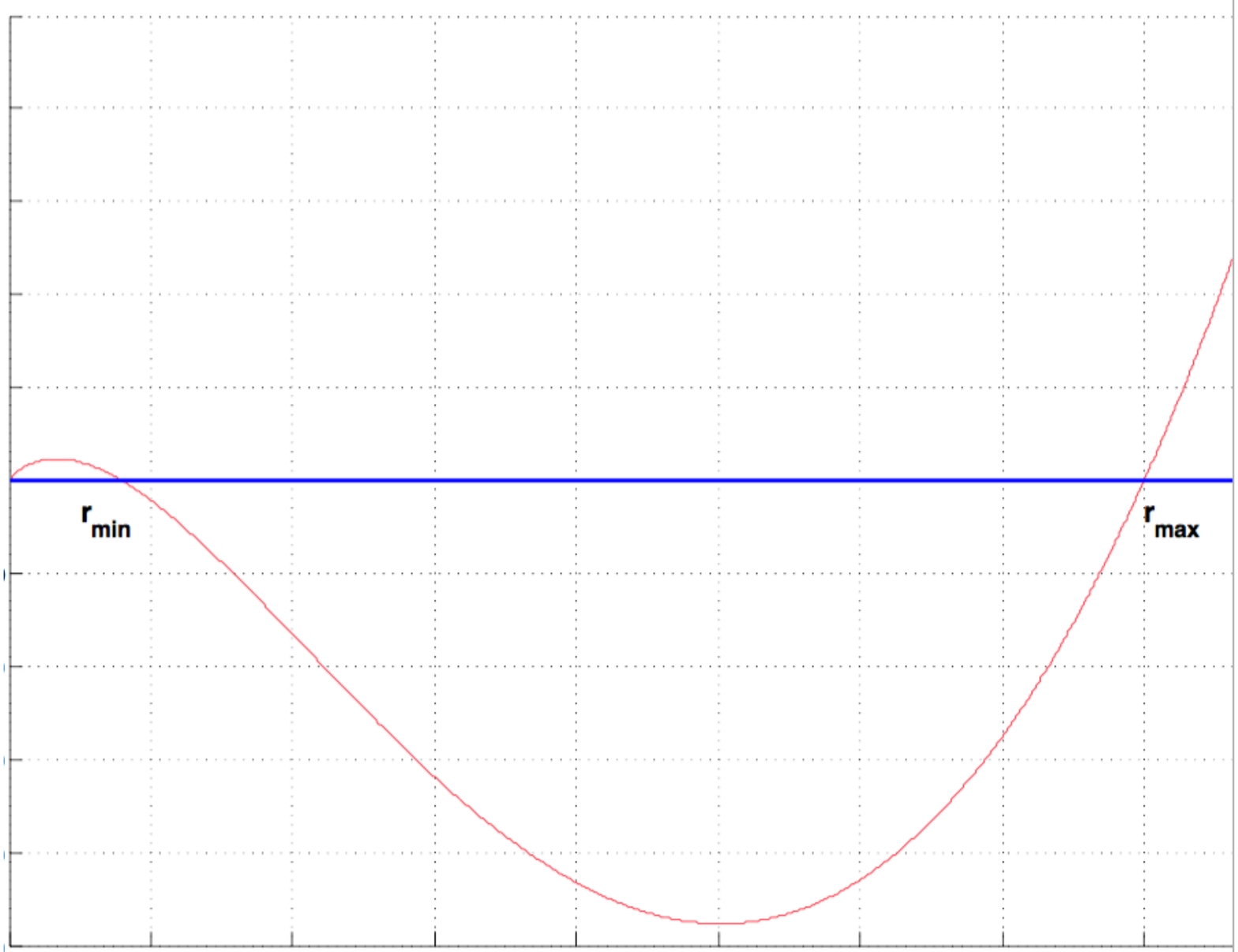}}
\caption{$r_{\tiny min}$ and $r_{\tiny max}$}
\end{figure}
If we also define $\Phi _{u}^{\prime \prime }(r)=\frac{\partial ^{2}}{%
\partial r^{2}}E_{\lambda }(ru)$, then, in case $\beta <1$ the equation $%
\Phi _{u}^{\prime }(r)=0$ may have at most two nonzero roots $r_{\min }>0$
and $r_{\max }>0$ such that $\Phi _{u}^{\prime \prime }(r_{\max })\geq 0$, $%
\Phi _{u}^{\prime \prime }(r_{\min })\leq 0$ and $0<r_{\max }\leq r_{\min }$
(see Figure 2), whereas, in case $\beta =1$ the equation $\Phi _{u}^{\prime
}(r)=0$ for any $\lambda >0$ has precisely one nonzero root $r_{\max }>0$
such that $\Phi _{u}^{\prime \prime }(r_{\max })\leq 0$. 
This implies that any weak solution of $SP(\alpha ,\beta ,\lambda )$ (any
critical point of $E_{\lambda }(u)$) corresponds to one of the cases $%
r_{\min }=1$ or $r_{\max }=1$. However, it was discovered in \cite{IlEg}
(see also \cite{ilDr}) that in case when we study compactly supported
solutions this correspondence essentially depends on the relation between $%
\alpha $, $\beta $ and $\rN$.

In the present paper, developing \cite{IlEg}, we introduce in the set of
relevant exponents $\mathcal{E}:=\{(\alpha ,\beta ):~0<\alpha <\beta\le 1\}$
the following critical exponents curve depending on the dimension $\rN$ 
\begin{equation}
\mathcal{C}(\rN):=\{(\alpha ,\beta )\in \mathcal{E}:~~2(1+\alpha )(1+\beta
)-\rN(1-\alpha )(1-\beta )=0\}.  \label{crN}
\end{equation}%
\begin{figure}[th]
\center{\includegraphics[width=0.8\linewidth]{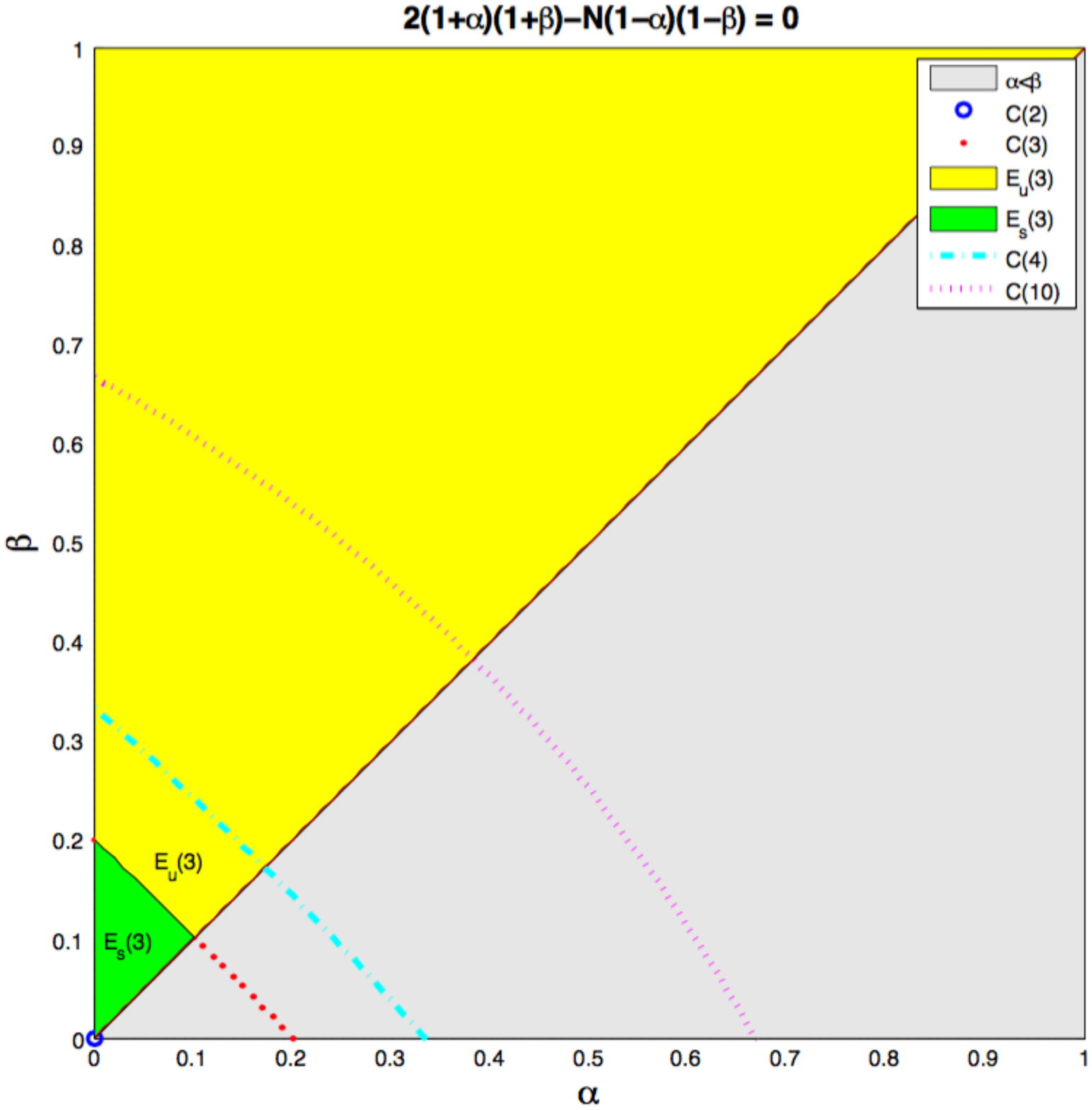}}
\caption{Sets $ \mathcal{E}_{s}(\rN)$ and  $\mathcal{E}_{u}(\rN)$ for $\rN=3,~4$ and 10}
\end{figure}
This curve exists if and only if $\rN\geq 3$ and it separates two sets of
exponents in $\mathcal{E}$ (see Figure 3) 
\begin{align*}
& \mathcal{E}_{s}(\rN):=\{(\alpha ,\beta )\in \mathcal{E}:~~2(1+\alpha
)(1+\beta )-\rN(1-\alpha )(1-\beta )<0\}, \\
& \\
& \mathcal{E}_{u}(\rN):=\{(\alpha ,\beta )\in \mathcal{E}:~~2(1+\alpha
)(1+\beta )-\rN(1-\alpha )(1-\beta )>0\},
\end{align*}%
whereas in the cases $\rN=1,2$ one has $\mathcal{E}=\mathcal{E}_{u}(\rN)$.

The main property of $\mathcal{C}(\rN)$ is contained in

\begin{lemma}
\label{lemG} Let $\rN\geq 1$ and let $\Omega $ be a bounded and star-shaped domain
in $\mathbb{R}^{\rN}$ whose boundary $\partial \Omega $ is a $\cC^{1}$-manifold.

\begin{description}
\item[1)] Assume $(\alpha ,\beta )\in \mathcal{C}(\rN)$.\thinspace\ ~ Then any
flat ground state solution $u$ of \eqref{1} satisfies $\Phi _{u}^{\prime \prime
}(r)|_{r=1}=0$.

\item[2)] Assume $(\alpha ,\beta )\in \mathcal{E}_{u}(\rN)$. Then any flat ground state
solution $u$ of \eqref{1} satisfies $\Phi _{u}^{\prime \prime }(r)|_{r=1}<0$.

\item[3)] Assume $(\alpha ,\beta )\in \mathcal{E}_{s}(\rN)$. Then any ground state
solution $u$ of \eqref{1} satisfies $\Phi _{u}^{\prime \prime }(r)|_{r=1}>0$.
\end{description}
\end{lemma}

\bigskip 

The existence of flat (or compactly supported) ground state solutions of (%
\ref{1}) in the case $\beta <1$, $\rN\geq 3$ and $(\alpha ,\beta )\in \mathcal{%
E}_{s}(\rN)$ has been obtained in \cite{IlEg}. Furthermore, the existence of
flat solutions of (\ref{1}) (not necessary ground states) in case $\rN\geq 1$, 
$0<\alpha <\beta \leq 1$ has been proved in \cite{diaz, diaz-H-Il,
Kaper1,Kaper2}. 

As already mentioned, one of the main goals of this paper is to study the $%
\rH_{0}^{1}$-stability of flat ground state solutions of $SP(\alpha ,\beta
,\lambda )$. We recall that, if $v(t;v_{0})$ is a weak solution
to $PP(\alpha ,\beta ,\lambda ,v_{0}),$ we shall say that $v(t;v_{0})$ is 
\textit{\ }$\rH_{0}^{1}$-\textit{stable} if, given any $\varepsilon >0$, there
exists $\delta >0$ such that 
\begin{equation}
||v(t;v_{0})-v(t;w_{0})||_{1}<\varepsilon ~~\mbox{for any}~w_{0}~\mbox{such that}%
~~||v_{0}-w_{0}||_{1}<\delta ,~~\forall t>0,
\end{equation}%
where we used the $\rH_{0}^{1}(\Omega )-$norm 
\begin{equation*}
||u||_{1}=\left (\int_{\Omega }|\nabla u|^{2}\,dx\right )^{1/2}.
\end{equation*}%
Conversely, we say that a solution $v(t;v_{0})$ of $PP(\alpha ,\beta
,\lambda ,v_{0})$ is $\rH_{0}^{1}$-\textit{unstable } if there is $\varepsilon
>0$ such that for any $\delta >0$ and $T>0$, there exists 
\begin{equation*}
w_{0}\in U_{\delta }(v_{0}):=\{w\in H_{0}^{1}(\Omega
):~||v_{0}-w||_{1}<\delta \}
\end{equation*}%
and there exists $T>0$ such that for any $t>T$ 
\begin{equation}
||v(t;v_{0})-v(t;w_{0})||_{1}>\varepsilon ,
\end{equation}%
where $v(t;w_{0})$ is any weak solution of $PP(\alpha ,\beta ,\lambda ,w_{0})
$. Furthermore, we will use also the following definition: a solution $%
u_{\lambda }$ of $SP(\alpha ,\beta ,\lambda )$ is said to be \textit{%
linearly unstable stationary solution} if $\lambda _{1}(-\Delta +\alpha
u_{\lambda }^{\alpha -1}-\lambda \beta u_{\lambda }^{\beta -1})<0$.

In what follows, we will also use the following definition (\cite{BR}, \cite%
{Hernandez-Mancebo-Vega}): a solution $v(t;v_{0})$ of $PP(\alpha ,\beta
,\lambda ,v_{0})$ is said to be \textit{globally }$\rH_{0}^{1}(\Omega )$-%
\textit{unstable }if for any $\delta >0$ there exists 
\begin{equation*}
w_{0}\in U_{\delta }(v_{0}):=\{w\in H_{0}^{1}(\Omega
):~||v_{0}-w||_{1}<\delta \}
\end{equation*}%
such that 
\begin{equation}
||v(t;v_{0})-v(t;w_{0})||_{1}\rightarrow \infty ~~\mbox{as}~~t\rightarrow
\infty .
\end{equation}

Motivated by the uniqueness results for the $PP(\alpha ,\beta ,\lambda
,w_{0})$, we shall assume later the following "isolation assumption":

\begin{description}
\item[\textbf{(U)}] Given $u_{\lambda }$ nonnegative ground state solution of $%
SP(\alpha ,\beta ,\lambda )$, there exists a "positive-neighborhood" 
\begin{equation*}
U_{\delta }(u_{\lambda }):=\{v\in H_{0}^{1}(\Omega ),\text{ }v\geq 0\text{
on }\Omega ,\text{ such that }~||u_{\lambda }-v||_{1}<\delta \},
\end{equation*}%
with $\delta >0$ such that $SP(\alpha ,\beta ,\lambda )$ has no other
non-negative weak solution in $U_{\delta }(u_{\lambda })\setminus u_{\lambda
}$. 
\end{description}

Our first two results concern the existence and (un-)stability of ground
states of (\ref{1}). In case $0<\alpha<\beta< 1$ we have

\begin{theo}
\label{Th1} Let $\rN\geq 1$, $0<\alpha <\beta <1$, $\Omega $ be a bounded
domain in $\mathbb{R}^{\rN}$, with a smooth boundary. Then

\begin{description}
\item[(1)] There exists $\lambda ^{\ast }>0$ such that for all $\lambda
>\lambda ^{\ast }$ problem (\ref{1}) has a ground state $u_{\lambda }$ which
is nonnegative in $\Omega $ and $u_{\lambda }\in \cC^{1,\kappa }(\overline{%
\Omega })\cap \cC^{2}(\Omega )$ for some $\kappa \in (0,1)$.

\item[(2)] Assume {\rm ({\bf U})}, then the ground state $u_{\lambda }$ is a $%
\rH_{0}^{1}(\Omega )$-stable stationary solution of the parabolic problem~\eqref{p1}.
\end{description}
\end{theo}

\bigskip 

In case $\beta=1$ we have

\begin{theo}
\label{Th1.1} Let $\rN\geq 1$, $\beta =1$, $0<\alpha <1$, $\Omega $ be a
bounded star-shaped domain in $\mathbb{R}^{\rN}$, with a smooth boundary. Then

\begin{description}
\item[(1)] There exists $\lambda^*>0$ such that for all $\lambda >\lambda^*$
problem (\ref{1}) has a ground state $u_\lambda$ which is nonnegative in $%
\Omega$ and $u\in \cC^{1,\kappa}(\overline{\Omega}) \cap \cC^2(\Omega)$ for some 
$\kappa \in(0,1)$.

\item[(2)] Assume {\rm ({\bf U})}, the ground state $u_{\lambda }$ is a globally $%
\rH_{0}^{1}(\Omega )$-unstable stationary solution of parabolic problem~\eqref{p1}.
\end{description}
\end{theo}

\bigskip 

Our main result on the $\rH_{0}^{1}(\Omega )$-stability and $\rH_{0}^{1}(\Omega )
$-unstability of flat ground state solutions for $0<\alpha <\beta ~<~1$ is the
following

\begin{theo}
\label{Th2} Let $\rN\geq 1$, $\Omega $ be a bounded domain in $\mathbb{R}^{\rN}$
whose boundary $\partial \Omega $ is a $\cC^{1}$-manifold.

\begin{description}
\item[(I)] Assume $\rN=1,2$. Then for every $(\alpha ,\beta )\in \mathcal{E}$
( i.e. $0<\alpha <\beta $) any flat ground state solution $u_{\lambda }$ of \eqref{1} is
a linearized unstable stationary solution of parabolic problem \eqref{p1}.

\item[(II)] Assume {\rm ({\bf U})}, $\rN\geq 3$ and $(\alpha ,\beta )\in \mathcal{E}_{u}(\rN)
$. Then any flat ground state solution $u_{\lambda }$ of \eqref{1} is a linearized
unstable stationary solution of the parabolic problem \eqref{p1}.

\item[(III)] Assume $\rN\geq 3$, $(\alpha ,\beta )\in \mathcal{E}_{s}(\rN)$ and $%
\Omega $ is a strictly star-shaped domain with respect to the origin. Then

\begin{description}
\item[(1)] there exists $\lambda ^{\ast }>0$ such that \eqref{1} has a flat
ground state $u_{\lambda ^{\ast }}$, $u_{\lambda ^{\ast }}\geq 0$ and $%
u_{\lambda ^{\ast }}\in \cC^{1,\gamma }(\overline{\Omega })\cap \cC^{2}(\Omega )$
for some $\gamma \in (0,1)$;

\item[(2)] If in addition {\rm ({\bf U})} holds then the flat ground state solution $%
u_{\lambda ^{\ast }}$ is a $\rH_{0}^{1}(\Omega )$-stable stationary solution
of the parabolic problem \eqref{p1}.
\end{description}
\end{description}
\end{theo}

\bigskip 

In the case $\beta=1$ we have

\begin{theo}
\label{Th2.2} Assume $\rN\geq 1$, $0<\alpha <1$, $\beta =1$ and $\Omega $ be a
bounded domain in $\mathbb{R}^{\rN}$ whose boundary $\partial \Omega $ is a $%
\cC^{1}$-manifold. Then

\begin{description}
\item[(1)] there exists $\lambda^*>0$ such that \eqref{1} has a ground state 
$u_{\lambda^*}$ which is a flat solution in $\Omega$ and $u_{\lambda^*}\geq
0 $ and $u_{\lambda^*}\in C^{1,\alpha}(\overline{\Omega}) \cap C^2(\Omega)$
for some $\alpha\in(0,1)$;

\item[(2)] If in addition {\rm ({\bf U})} holds the flat ground state solution $%
u_{\lambda ^{\ast }}$ is globally $\rH_{0}^{1}(\Omega )$-unstable stationary
solution of the parabolic problem \eqref{p1}.
\end{description}
\end{theo}

\bigskip 

The limit case $\alpha =0$ can be also considered. In particular, this shows that the first \textquotedblleft compressed
mode\textquotedblright\ function (solution of $SP(0,1,\lambda )$: see \cite%
{Ozolis-Osher}, \cite{Ozolis-Osher 2}), of great relevance in signal
processing,  is globally $\rH_{0}^{1}(\Omega )$-unstable.

\section{Parabolic problem. Existence, uniqueness and boundness on
non-negative solutions}

Given $v_{0}\in \rL^{\infty }(\Omega ),$ $v_{0}\geq 0,$ we shall say that $%
v\in \cC([0,+\infty ),\rL^{2}(\Omega ))$ is a weak solution of $PP(\alpha ,\beta
,\lambda ,v_{0})$ if $v\geq 0$, $\lambda v^{\beta }-v^{\alpha }\in \rL^{\infty
}((0,T)\times \Omega ),$ for any $T>0$ and 
\begin{equation}
v(t)=T(t)v_{0}+\int_{0}^{t}T(t-s)(\lambda v^{\beta }(s)-v^{\alpha }(s))ds.
\end{equation}%
Here $(T(t))_{t\geq 0}$ is the heat semigroup with homogeneous Dirichlet
boundary conditions, i.e. $T(t)=e^{t(-\Delta )}$. The existence of weak
solutions is an easy variation of previous results in the literature (see,
e.g. \cite{CDE}, \cite{Akagi} and the works \cite{Davila-Montenegro survey},
\cite{AnhDiazPaul} dealing with the more difficult case of singular
equations $\alpha \in (-1,0)$). For the reader convenience we shall collect
here some additional regularity information on weak solutions of $%
PP(\alpha ,\beta ,\lambda ,v_{0}).$

\bigskip

\begin{prop} \label{proposition2.1}
For any $v_{0}\in \rL^{\infty }(\Omega ),$ $v_{0}\geq 0$ there exists a
nonnegative weak solution $v\in \cC([0,+\infty ),\rL^{2}(\Omega ))$ of $%
PP(\alpha ,\beta ,\lambda ,v_{0})$. In fact, for every $p\in \lbrack
1,\infty ],$ $v\in \cC([0,+\infty );\rL^{p}(\Omega ))$, and if $p<\infty $ 
\begin{equation}
v-T(.)v_{0}\in \rL^{p}(\tau ,T;\rW^{2,p}(\Omega )\cap \rW_{0}^{1,p}(\Omega ))\cap
\rW^{1,p}(\tau ,T;\rL^{p}(\Omega )),  \label{additional regularity}
\end{equation}%
for any $0<\tau <T.$ In particular, $v$ satisfies the equation $PP(\alpha
,\beta ,\lambda ,v_{0})$ for a.e. $t\in (0,+\infty )$. Moreover, if we also
assume that $v_{0}\in \rH_{0}^{1}(\Omega )$ then $\frac{\partial }{\partial t}%
E_{\lambda }(v(.))\in \rL^{1}(\tau ,T)$, function $E_{\lambda }(v(.))$ is
absolutely continuous for a.e. $t\in (\tau ,T$) 
\begin{equation}
\frac{\partial }{\partial t}E_{\lambda }(v(t))=\int_{\Omega }\big (\lambda
v^{\beta}+v|^{\alpha})v_{t}(t)dx-\int_{\Omega }v_{t}(t)^{2}dx.
\label{derivada de E}
\end{equation}
\end{prop}

\proof Among many possible methods\textit{\ }to prove the existence
of weak solutions we shall follow here the one based on a fixed point
argument as in \cite{Diaz-Vrabie} (see also \cite{Diaz-Tello} where the case 
$\beta =0$ was considered on a Riemannian manifold). For every $h\in
\rL^{\infty }((0,T)\times \Omega )$ we consider the problem $(P_{h})$

\begin{equation*}
(P_{h})\quad \left\{ 
\begin{array}{ll}
v_{t}-\Delta v+|v|^{\alpha -1}v=h & \text{in }(0,+\infty )\times \Omega  \\ 
v=0 & \text{on }(0,+\infty )\times \partial \Omega  \\ 
v(0,x)=v_{0}(x) & \text{on }\Omega ,%
\end{array}%
\right. 
\end{equation*}%
which we can reformulate in terms of an abstract Cauchy problem on the
Hilbert space $\rH=\rL^{2}(\Omega )$ as

\begin{equation*}
(P_{h})=\left\{ 
\begin{array}{lr}
\dfrac{dv}{dt}(t)+\mathcal{A}v(t)=h(t) & t\in (0,T)\text{, in }\rH, \\ [.2cm]
v(0)=v_{0} & 
\end{array}%
\right. 
\end{equation*}%
where $\mathcal{A}=\partial \varphi $ denotes the subdifferential of the
convex function%
\begin{equation*}
\varphi (v)=\left\{ 
\begin{array}{ll}
\displaystyle \frac{1}{2}\int_{\Omega }|\nabla v|^{2}\,dx+\frac{1}{\alpha +1}\int_{\Omega
}|v|^{\alpha +1}\,dx & \text{if }v\in \rH_{0}^{1}(\Omega )\cap \rL^{\alpha
+1}(\Omega ) \\ 
+\infty  & \text{otherwise,}%
\end{array}%
\right. 
\end{equation*}%
(see, e.g. \cite{brrr}, \cite{brezis1971} and \cite{Diaz-vol-1}). As in \cite%
{Diaz-Vrabie}, \cite{Diaz-Tello}\textbf{, }we define the operator $\mathcal{T%
}:h\rightarrow g$ where $g=\lambda |v_{h}|^{\beta -1}v_{h}$ and $v_{h}$ is
the solution of $(P_{h})$. It is easy to see that every fixed point of $%
\mathcal{T}$ is a solution of $PP(\alpha ,\beta ,\lambda ,v_{0}).$ Then $%
\mathcal{T}$ satisfies the hypotheses of Kakutani Fixed Point Theorem (see
e.g. Vrabie \cite{vrabie}), since if $X=\rL^{2}((0,T),\rL^{2}(\Omega
))$ then

\begin{itemize}
\item[(i)] $K=\{h\in \rL^{2}\big (0,T,\rL^{\infty }(\Omega )\big ):||h(t)||_{\rL^{\infty}(\Omega)}\leq C_{0}$
a.e. $t\in (0,T)\}$ is a nonempty, convex and weakly compact set of $X;$

\item[(ii)] $\mathcal{T}:K\mapsto 2^{X}$ with nonempty, convex and closed
values such that $\mathcal{T}(g)\subset K$, $\forall g\in K;$

\item[(iii)] graph($\mathcal{T}$) is weakly$\times $weakly sequentially
closed.
\end{itemize}

\noindent Consequently, $\mathcal{T}$ \ has at least one fixed point in $K$
which is a local (in time) solution of $PP(\alpha ,\beta ,\lambda ,v_{0})$.
The final key point is to show that there is no blow-up phenomenon. This hods by the a
priori estimate 
\begin{equation*}
0\leq v(t,x)\leq z(t,x)\text{, for any }t\in \lbrack 0,+\infty )\times
\Omega ,
\end{equation*}%
where $v(t,x)$ is any weak solution of $PP(\alpha ,\beta ,\lambda ,v_{0})$
and $z(t,x)$ is the solution of the corresponding auxiliary problem 
\begin{equation}
\left\{ 
\begin{array}{ll}
z_{t}-\Delta z=\lambda z^{\beta } & \text{in }(0,+\infty )\times \Omega  \\ 
z=0 & \text{on }(0,+\infty )\times \partial \Omega  \\ 
z(0,x)=v_{0}(x) & \text{on }\Omega .%
\end{array}%
\right. 
\end{equation}%
This implies that there is no finite blow-up (and thus the maximal existence time is $%
T_{\max }=+\infty $). In particular, if $\beta \in (0,1)$ we have the
estimate

\begin{equation*}
\left\Vert v(t)\right\Vert _{\rL^{\infty }(\Omega )}\leq \big (\left\Vert
v_{0}\right\Vert _{\rL^{\infty }(\Omega )}^{1-\beta }+(1-\beta )t\big )^{1/(1-\beta
)}.
\end{equation*}%
If $\beta =1$ then the function $w(t,x)=v(t,x)e^{-\lambda t}$ satisfies%
\begin{equation}
\left\{ 
\begin{array}{ll}
w_{t}-\Delta w+e^{-\lambda (1-\alpha )t}w^{\alpha }=0 & \text{in }(0,+\infty
)\times \Omega \\ 
w=0 & \text{on }(0,+\infty )\times \partial \Omega \\ 
w(0,x)=v_{0}(x) & \text{on }\Omega ,%
\end{array}%
\right.  \label{Extincion con peso}
\end{equation}%
which is uniformly (pointwise) bounded by the solution of the linear heat
equation with the same initial datum. Since the operator $A=\overline{%
\partial \varphi }^{\rL^{p}(\Omega )\times \rL^{p}(\Omega )}$ is m-accretive in $%
\rL^{p}(\Omega )$ for every $p\in \lbrack 1,\infty ]$ (see, e.g. the
presentation made in \cite{Diaz-vol-1}), by the regularity results for
semilinear accretive operators we conclude the first part of the additional
regularity of the statement (\ref{additional regularity}). Finally, by
Theoreme 3.6 of \cite{brrr} we know that $\frac{\partial }{\partial t}%
\varphi (v_{h})\in \rL^{1}(\tau ,T),~\varphi
(v_{h})$ is absolutely continuous and for a.e. $t\in (\tau ,T$) 
\begin{equation*}
\frac{\partial }{\partial t}\varphi (v_{h})=\int_{\Omega
}(h(t))(v_{h})_{t}(t)dx-\int_{\Omega }\left[ (v_{h})_{t}(t)\right] ^{2}dx.
\end{equation*}%
Then (\ref{derivada de E}) holds by taking $h=\lambda |v_{h}|^{\beta
-1}v_{h}$ (the fixed point of $\mathcal{T}$). \fin
\newline

\bigskip

\begin{coro}
Assume $\beta =1$. Then the weak solution is unique.
\end{coro}

\proof Thanks to the change of variable $%
w(t,x)=v(t,x)e^{-\lambda t}$ the problem becomes (\ref{Extincion con peso})
and the result follows from the semigroup theory since it is well-known
(see, e.g., \cite{diaz} Chapter 4) that the operator $Aw:=-\Delta
w+e^{-\lambda (1-\alpha )t}\left\vert w\right\vert ^{\alpha -1}w$ is a
T-accretive operator in $\rL^{p}(\Omega )$ for any $p\in \lbrack 1,+\infty ].$ \fin

\bigskip

A more delicate question deals with the proof of the uniqueness of weak
solutions for $\beta \in (0,1)$. We point out that some previous results in
the literature dealing with the case $\beta \in (0,1)$ (see \cite{CDE} and
its references) are not applicable to our framework due to the presence of
the absorption term $|v|^{\alpha -1}v$.

\bigskip

We define the following class of functions: 
\begin{multline}
\mathcal{M}(\nu ,T):=\Big\lbrace v\in \mathcal{C}\left( [0,T];\rL^{2}(\Omega
)\right) \;\big\vert\;\forall \,T^{\prime }\in
(0,T),\;\hbox{there exists $C(T^{\prime })>0$ such that:}  \label{M(nu)} \\
\forall t\in (0,T^{\prime }),\;v(t,x)\geq C(T^{\prime })d(x)^{\nu }\quad 
\text{in }\Omega \Big\rbrace,
\end{multline}%
where $\delta (x):=\mathrm{dist}(x,\partial \Omega )$ (which we shall denote
simply as $\delta $) and 
\begin{equation}
\nu \in \left( 0,\frac{2}{1-\alpha }\right] .
\end{equation}

The following result collects some useful estimates leading to the
uniqueness of non-degenerate weak solutions: 

\bigskip

\begin{theo}
Let $w$ (resp. $v$) be a weak subsolution $PP(\alpha ,\beta ,\lambda ,w_{0})$%
, i.e. 
\begin{equation*}
\left\{ 
\begin{array}{ll}
w_{t}-\Delta w+|w|^{\alpha -1}w\leq \lambda |w|^{\beta -1}w & \text{in }%
(0,+\infty )\times \Omega  \\ 
w=0 & \text{on }(0,+\infty )\times \partial \Omega  \\ 
w(0,x)=w_{0}(x) & \text{on }\Omega ,%
\end{array}%
\right. 
\end{equation*}%
with $w\in \mathcal{C}\left( [0,T];\rL^{2}(\Omega )\right) \cap \rL^{\infty
}((0,T)\times \Omega )\cap \rL_{loc}^{2}(0,T:\rH_{0}^{1}(\Omega ))$, $w\in
\rH_{loc}^{1}(0,T:\rH^{-1}(\Omega ))$ (resp. similar conditions for $v$ but
with the reversed inequalities). Then:

\noindent i) If $v\in \mathcal{M}(\nu ,T)$ for some $\nu \in \left( 0,\frac{2%
}{1-\alpha }\right] $, there exists a constant $C>0$ such that for any $t\in
\lbrack 0,T)$, we have 
\begin{equation}
\left\Vert \lbrack w(t)-v(t)]_{+}\right\Vert _{{\rL}^{2}(\Omega )}\leq
e^{\lambda Ct}\left\Vert [w_{0}-v_{0}]_{+}\right\Vert _{\rL^{2}(\Omega )}.
\label{CDE}
\end{equation}

\noindent ii) If $w\in \mathcal{M}(\nu ,T)$ for some $\nu \in \left( 0,\frac{%
2}{1-\alpha }\right] ,$ there exists a constant $C>0$ such that for any $%
t\in \lbrack 0,T)$, we have 
\begin{equation}
\left\Vert \lbrack w(t)-v(t)]_{-}\right\Vert _{{L}^{2}(\Omega )}\leq
e^{\lambda Ct}\left\Vert [w_{0}-v_{0}]_{-}\right\Vert _{L^{2}(\Omega )}.
\end{equation}

\noindent iii) Assume $w_{0}\leq v_{0}$ and that $v\in \mathcal{M}(\nu ,T)$
or $w\in \mathcal{M}(\nu ,T).$ Then, for any $t\in \lbrack 0,T]$, $w(t,\cdot
)\leq v(t,\cdot )\quad $ a.e. in $\Omega .$

\noindent iv) There is uniqueness of weak solutions in the class $\mathcal{M}%
(\nu ,T)$. Moreover, if $v,w\in \mathcal{M}(\nu ,T)$ are weak solutions of $%
PP(\alpha ,\beta ,\lambda ,w_{0})$ and $PP(\alpha ,\beta ,\lambda ,v_{0})$,
respectively, then there exists a constant $C>0$ such that for any $t\in
\lbrack 0,T)$, we have 
\begin{equation}
\left\Vert w(t)-v(t)\right\Vert _{{\rL}^{2}(\Omega )}\leq e^{\lambda
Ct}\left\Vert w_{0}-v_{0}\right\Vert _{L^{2}(\Omega )}.
\end{equation}
\end{theo}

~\newline

We shall get later some sufficient conditions on the initial datum $v_{0}$
ensuring that there exists some weak solution of $PP(\alpha ,\beta ,\lambda
,v_{0})$ belonging to the class $\mathcal{M}(\nu ,T)$.

\bigskip
\noindent 
{\sc Proof of Theorem 2.1.}\quad Multiplying by $(w(t)-v(t))_{+}$
the difference of the inequalities satisfied by $w$ and $v$ we obtain 
$$
\begin{array}{c}
\hspace*{-2cm}\displaystyle \frac{1}{2}\frac{d}{dt}\int_{\Omega }[w(t)-v(t)]_{+}^{2}+\int_{\Omega
}\left\vert \nabla \lbrack w(t)-v(t)]_{+}\right\vert ^{2}+\int_{\Omega
}(w(t)^{\alpha }-v(t)^{\alpha })[w(t)-v(t)]_{+} \\ [.2cm]
\hspace*{8,5cm}\displaystyle \leq \lambda \int_{\{w>v\}}(w(t)^{\beta }-v(t)^{\beta })[w(t)-v(t)].
\end{array}
$$
But, since $\beta \in (0,1)$ 
\begin{equation*}
w^{\beta }-v^{\beta }\leq \frac{\beta }{v^{1-\beta }}(w-v)\text{ for any }%
0<v<w\leq M
\end{equation*}%
for some $M>0.$ On the other hand, since $v\in \mathcal{M}(\nu ,T),$ and $%
\alpha <\beta ,$ by applying Young's inequality we get 
\begin{equation*}
v^{\beta -1}\leq \frac{1}{C^{(1-\beta )}d(x)^{\nu (1-\beta )}}\leq \frac{%
\varepsilon }{d(x)^{2}}+C_{\varepsilon },
\end{equation*}%
for any $\varepsilon >0$ and for some $C_{\varepsilon }>0.$ Then, from the
monotonicity of the function $w\rightarrow $ $w^{\alpha }$, taking $M=\max
(\left\Vert w\right\Vert _{\rL^{\infty }((0,T)\times \Omega )},\left\Vert
v\right\Vert _{\rL^{\infty }((0,T)\times \Omega )})$ we obtain
$$
\frac{1}{2}\frac{d}{dt}\int_{\Omega }[w(t)-v(t)]_{+}^{2}+\int_{\Omega
}\left\vert \nabla \lbrack w(t)-v(t)]_{+}\right\vert ^{2} \leq \lambda \varepsilon \int_{\Omega }\frac{[w(t)-v(t)]_{+}^{2}}{d(x)^{2}}%
+\lambda C_{\varepsilon }\int_{\Omega }[w(t)-v(t)]_{+}^{2}.
$$
Applying Hardy's inequality,%
\begin{equation*}
\int_{\Omega }\frac{z^{2}}{d(x)^{2}}dx\leq C\int_{\Omega }\left\vert \nabla
z\right\vert ^{2}dx
\end{equation*}%
for any $z\in H_{0}^{1}(\Omega )$, choosing $\varepsilon >0$ sufficiently
small and using Gronwall's inequality we get the conclusion~i). The proof of
ii) is similar but this time we multiply by $(v(t)-w(t))_{-}$ the difference
of the inequalities satisfied by $v$ and $w$ and use the fact that, since $%
\beta \in (0,1)$, 
\begin{equation*}
(w^{\beta }-v^{\beta })[w(t)-v(t)]_{-}\leq \frac{\beta }{w^{1-\beta }}%
[w(t)-v(t)]_{-}^{2}\text{ for any }0<v,w\leq M,
\end{equation*}%
for some $M>0.$ Again, since $v\in \mathcal{M}(\nu ,T),$ and $\alpha <\beta ,
$ by applying Young's inequality we get 
\begin{equation*}
w^{\beta -1}\leq \frac{1}{C^{(1-\beta )}d(x)^{\nu (1-\beta )}}\leq \frac{%
\varepsilon }{d(x)^{2}}+C_{\varepsilon },
\end{equation*}%
for any $\varepsilon >0$ and for some $C_{\varepsilon }>0$ and the proof
ends as in the case i). The proofs of iii) and iv) are easy consequences of
i) and ii). \fin

\bigskip

\begin{prop}
\label{Proposition non degen}Assume 
\begin{equation}
v_{0}(x)\geq K_{0}d(x)^{2/(1-\alpha )}%
\hbox{for any
$x\in \overline{\Omega}$},  \label{initial datum}
\end{equation}%
for some constant $K_{0}>0.$ Let $v$ be a weak solution of $PP(\alpha ,\beta
,\lambda ,v_{0})$. Then:
\par
\noindent 
a) Given $T>0$ for any $K_{0}>0$ there is a $T_{0}=T_{0}(K_{0})\in (0,T]$
such that $v\in \mathcal{M}(\nu ,T_{0})$ on for $\nu =2/(1-\alpha ).$
\par
\noindent 
b) If $K_{0}$ and $\lambda $ are large enough then $v\in \mathcal{M}(\nu ,T)$
for $\nu =2/(1-\alpha )$, for any $T>0.$
\end{prop}
\proof By iii) of the above theorem it is enough to construct a
(local) subsolution satisfying the required boundary behavior.\textit{\ }We
shall carry out such construction by adapting the techniques presented in 
\cite{D ambiguity 2} (see also some related local subsolutions in \cite{AD}, 
\cite{DRM} and \cite{D ambiguity}). From the assumption (\ref{initial datum}%
) for any $x_{0}\in \partial \Omega $ there exist $\epsilon >0,\delta \geq 1,
$ $C_{0}>0$ and $x_{1}\in \Omega $ with $B_{\delta \epsilon }(x_{1})\subset
\Omega $ such that%
\begin{equation}
v_{0}(x)\geq C_{0}|x-x_{0}|^{\nu }\quad \text{\textit{a.e.}}\;x\in
B_{\epsilon }(x_{1}).  \label{Creu}
\end{equation}%
Let us take $x_{1}\in \Omega $ such that $\delta \epsilon >|x_{1}-x_{0}|\geq
((\delta +1)/2)\epsilon $, and define 
\begin{equation*}
\underline{U}(x)=\left\{ 
\begin{array}{ll}
K_{1}\epsilon ^{\nu }-K_{2}|x-x_{1}|^{\nu } & \text{\textit{if} }\;0\leq
|x-x_{1}|\leq \epsilon , \\ [.175cm]
K_{3}(\delta \epsilon -|x-x_{1}|)^{\nu } & \text{\textit{if} }\;\epsilon
\leq |x-x_{1}|\leq \delta \epsilon ,%
\end{array}%
\right. 
\end{equation*}%
and, for $x\in B_{\delta \epsilon }(x_{1})$ and $t\in (0,T]$

\begin{equation*}
\underline{V}(t,x)=\varphi (t)\underline{U}(x).
\end{equation*}%
We shall show that it is possible to choose all the above constants and
function $\varphi (t)$ such that $\underline{V}$ is a weak subsolution of $%
PP(\alpha ,\beta ,\lambda ,v_{0})$ with the desired growth near $\partial
B_{\delta \epsilon }(x_{1})$ for suitable time interval $[0,T_{0}(K_{0}))$
in case a) or on the whole interval $[0,T]$ in case b). Since $\underline{U}%
(x)=\eta (|x-x_{1}|)$ on $B_{\delta \epsilon }(x_{1})$ then the Laplacian
operator can be written as

$$
\Delta \eta (r)=\eta ^{\prime \prime }(r)+\frac{\rN-1}{r}\eta ^{\prime }(r)
$$
with $r\in (0,\delta \epsilon )$. By defining $\eta _{1}(r)=K_{1}\epsilon
^{\nu }-K_{2}r^{\nu }$ and $\eta _{2}(r)=K_{3}(\delta \epsilon -r)^{\nu }$
then 
\begin{equation*}
\eta (r)=\left\{ 
\begin{array}{ll}
\eta _{1}(r) & 0\leq r\leq \epsilon , \\ [.175cm]
\eta _{2}(r) & \epsilon \leq r\leq \delta \epsilon .%
\end{array}%
\right.
\end{equation*}%
The list of conditions which we must check to ensure that $\underline{V}%
(t,x) $ is a local-weak- subsolution is the following:

\noindent 1) $\underline{V}\in \mathcal{C}\left( [0,T];\rL^{2}(B_{\delta
\epsilon }(x_{1}))\right) \cap \rL^{\infty }((0,T)\times B_{\delta \epsilon
}(x_{1}))\cap \rL_{loc}^{2}(0,T:\rH_{0}^{1}(B_{\delta \epsilon }(x_{1})))$, $%
\underline{V}\in \rH_{loc}^{1}(0,T:\rH^{-1}(B_{\delta \epsilon }(x_{1}))).$
This is guarantied if we take $\varphi \in \rH^{1}(0,T)$ and $\underline{U}\in
\cC^{1}(B_{\delta \epsilon }(x_{1}))$ (since by construction $\underline{U}=0$
on $\partial B_{\delta \epsilon }(x_{1}).$ In particular, we must have%
\begin{equation}
(K_{1}-K_{2})\epsilon ^{\nu }=K_{3}(\epsilon (\delta -1))^{\nu }
\label{continuity}
\end{equation}%
\begin{equation}
\nu K_{2}\epsilon ^{\nu -1}=-\nu K_{3}(\epsilon (\delta -1))^{\nu -1}.
\label{C^1}
\end{equation}

\noindent 2) $\underline{V}(0,x)\leq v_{0}(x)$ a.e. on $B_{\delta \epsilon
}(x_{1}).$ Thanks to (\ref{Creu}), since $\eta _{1}(r)$ is concave and $%
C_{0}r^{\nu }$ is convex it is enough to have%
\begin{equation*}
\varphi (0)K_{3}(\epsilon (\delta -1))^{\nu }\leq C_{0}\epsilon ^{\nu }\text{
on }B_{\delta \epsilon }(x_{1}).
\end{equation*}

\noindent 3) $\underline{V}_{t}-\Delta \underline{V}+\underline{V}^{\alpha
}\leq \lambda \underline{V}^{\beta }$ (in a weak form) on $%
[0,T_{0}(K_{0}))\times B_{\delta \epsilon }(x_{1}).$ For $\mu >0$ let us
introduce $\mathcal{L}(\eta :\mu )=-\Delta \eta +\mu \eta ^{\alpha }$. Then,
if we write $r=\epsilon s$ 
\begin{equation*}
\begin{array}{l}
\mathcal{L}(\eta _{1})\leq \nu (\nu -1)K_{2}r^{\nu -2}+\nu (\rN-1)K_{2}r^{\nu
-2}+\mu \lbrack K_{1}\epsilon ^{\nu }-K_{2}r^{\nu }]^{\alpha }\quad \\ [.2cm]
=\left[ \nu (\nu -1)K_{2}s^{\nu \alpha }+\nu (\rN-1)K_{2}s^{\nu \alpha }+\mu
(K_{1}-K_{2}s^{\nu })^{\alpha }\right] \epsilon ^{\alpha \nu }\\ [.2cm]
\leq K_{4}\epsilon ^{\nu \alpha }%
\end{array}%
\end{equation*}%
where 
\begin{equation}
K_{4}=K_{4}(\mu ):=\nu \lbrack (\nu -1)+(\rN-1)K_{2}]+\mu K_{1}.  \label{K_4}
\end{equation}

\noindent On the other hand, 
\begin{equation*}
\begin{array}{l}
\displaystyle \mathcal{L}(\eta _{2})\leq -\lambda \nu (\nu -1)K_{3}(\delta \epsilon
-r)^{\nu -2}+(\rN-1)\nu K_{3}\frac{(\delta \epsilon -r)^{\nu -1}}{r}+\mu
K_{3}^{\alpha }(\delta \epsilon )^{\nu \alpha } \\ [.2cm]
\displaystyle \leq \nu K_{3}(\delta \epsilon -r)^{\nu \alpha }\left( -(\nu -1)+(\rN-1)\frac{%
(\delta \epsilon -r)}{r}+\mu K_{3}^{\alpha -1}\nu ^{-1}\right) .%
\end{array}%
\end{equation*}%
Now $\dfrac{(\delta \epsilon -r)}{r}\leq \delta -1$ when $\epsilon \leq r\leq
\delta \epsilon $ and thus if 
\begin{equation}
1\leq \delta <1+(\nu \alpha +1)/(\rN-1)  \label{delta}
\end{equation}

\noindent so, if we choose $K_{3}$ as 
\begin{equation}
K_{3}=K_{3}(\mu ,\delta ):=\left( \frac{\mu \nu ^{-1}}{(\nu \alpha
+1)-(\rN-1)(\delta -1)}\right) ^{\frac{1}{1-\alpha }},  \label{K-3}
\end{equation}%
we obtain that $-\Lambda \eta _{2}+\mu \eta _{2}^{\alpha }\leq 0$.

\noindent Moreover,%
\begin{equation*}
\underline{V}_{t}-\Delta \underline{V}+\underline{V}^{\alpha }=\varphi
^{\prime }\eta -\varphi \left (\eta ^{\prime \prime }+\frac{\rN-1}{r}\eta ^{\prime
}\right )+\varphi ^{\alpha }\eta ^{\alpha }.
\end{equation*}%
Then, if we have $\varphi \in \cC^{1}(0,T)$ such that 
\begin{equation}
\varphi ^{\prime }(t)\leq 0,  \label{fi decreasing}
\end{equation}%
then once we have 
\begin{equation}
\varphi (0)\leq 1,  \label{fi en cero}
\end{equation}%
given $\varepsilon _{1}\in (0,1)$, we always can find $T_{0}(\varepsilon
_{1})\leq T$ such that 
\begin{equation*}
\varepsilon _{1}\leq \varphi (t)\leq 1\text{ for any }t\in \lbrack
0,T_{0}(\varepsilon _{1})]
\end{equation*}%
and hence, if 
\begin{equation}
\mu =\frac{1}{(\varepsilon _{1})^{1-\alpha }}  \label{mu}
\end{equation}

\begin{equation*}
\Delta \underline{V}+\underline{V}^{\alpha }\leq (\varepsilon
_{1})^{1-\alpha }\varphi (t)^{\alpha }(-\Delta \eta (r)+\mu \eta ^{\alpha
})\leq 0.
\end{equation*}%
This implies that $\underline{V}_{t}-\Delta \underline{V}+\underline{V}%
^{\alpha }\leq \lambda \underline{V}^{\beta }$ (in a weak form) on $%
[0,T_{0}(\varepsilon _{1}))\times (B_{\delta \epsilon }(x_{1})\setminus
B_{\epsilon }(x_{1}))$. The remaining condition is to have the above
inequality also on $B_{\epsilon }(x_{1})$. This will be an easy consequence
if we take as function $\varphi $ any subsolution of the associated ODE:
more precisely. such that 
\begin{equation*}
\varphi ^{\prime }(t)+\frac{(\max \eta _{1})^{\alpha }}{\min \eta _{1}}%
\varphi (t)^{\alpha }\leq \frac{\lambda }{(\min \eta _{1})^{1-\beta }}%
\varphi (t)^{\beta }.
\end{equation*}%
By taking $\varphi (0)$ and $\varepsilon _{1}$ small enough it is easy to
see that it is possible to choose the rest of constants such that all the
above conditions follow and this ends the proof of case a). In case b) the
arguments are very similar but in this case it is possible to take as
function $\varphi (t)$ the one given by 
\begin{equation*}
\varphi (t)=(\varepsilon _{2}+e^{-kt})
\end{equation*}%
for suitable $\varepsilon _{2}>0$ and $k>0$ small enough. \fin

\bigskip 

\begin{coro}
Assume $v_{0}$ as in Proposition 2.2 and let $v$ be a weak solution of $%
PP(\alpha ,\beta ,\lambda ,v_{0})$ such that the nondegeneracy constant $C$
in (\ref{M(nu)}) is independent of $T$, for any $T>0.$ Let $u\in \rL^{\infty
}(\Omega )$ be a solution of the stationary problem $SP(\alpha ,\beta
,\lambda )$ such that $v(t)\rightarrow u$ in $\rL^{2}(\Omega )$ a.e. $%
t\nearrow +\infty .$ Then $u$ satisfies the nondeneracy property $u(x)\geq
Kd(x)^{2/(1-\alpha )}$ for some $K>0.$
\end{coro}

\bigskip 

The stability of the trivial solution $u\equiv 0$ of $SP(\alpha ,\beta
,\lambda )$ for $\lambda $ small is very well illustrated by means of the
following "extinction in finite time" property of solutions of the
associated parabolic problem $PP(\alpha ,\beta ,\lambda ,v_{0})$ assumed $%
\lambda $ small enough.

\bigskip

\begin{theo}
Assume 
\begin{equation}
\lambda \in \lbrack 0,\lambda _{1}).  \label{hypo extinction}
\end{equation}%
Let $v_{0}\in \rL^{\infty }(\Omega ),$ $v_{0}\geq 0.$ Assume $\beta =1$ or (%
\ref{initial datum}). Then there exists $T_{0}>0$ such that the solution $v$
of $PP(\alpha ,\beta ,\lambda ,v_{0})$ satisfies that $v(t)\equiv 0$ on $%
\Omega $ for any $t\geq T_{0}.$
\end{theo}

\proof We shall use an energy method in the spirit of 
\cite{ADS} (see also \cite{GiSaSer}). By multiplying by $v(t)$ and
integrating by parts (as in the proof of uniqueness) we arrive to%
\begin{equation*}
\frac{1}{2}\frac{d}{dt}\int_{\Omega }v(t)^{2}dx+\int_{\Omega }\left\vert
\nabla v(t)\right\vert ^{2}dx+\int_{\Omega }v(t)^{\alpha +1}dx=\lambda
\int_{\Omega }v(t)^{\beta +1}dx.
\end{equation*}%
Assume now that $\beta =1$. Then, by using the Poincar\'{e} inequality 
\begin{equation}
\lambda _{1}\int_{\Omega }v(t)^{2}dx\leq \int_{\Omega }\left\vert \nabla
v(t)\right\vert ^{2}dx  \label{Poincare}
\end{equation}%
we get%
\begin{equation*}
\frac{1}{2}\frac{d}{dt}\int_{\Omega }v(t)^{2}dx+\left (1-\frac{\lambda }{\lambda
_{1}}\right )\int_{\Omega }\left\vert \nabla v(t)\right\vert ^{2}dx+\int_{\Omega
}v(t)^{\alpha +1}dx\leq 0
\end{equation*}%
and the result holds exactly as in Proposition 1.1, Chapter 2 of \cite{ADS}. Indeed,
by applying the Gagliardo-Nirenberg inequality,%
\begin{equation*}
\left[ \int_{\Omega }v^{r}dx\right] ^{1/r}\leq C\left[ \int_{\Omega
}\left\vert \nabla v\right\vert ^{2}dx\right] ^{\theta /2}\left[
\int_{\Omega }vdx\right] 
\end{equation*}%
for any $r\in \lbrack 1,+\infty )$ if $\rN\leq 2$ and $r\in \left [ 1,\dfrac{2\rN%
}{\rN-2}\right ]$ if $\rN>2$ (with $\theta =\dfrac{2\rN(r-1)}{ r+2\rN} \in (0,1)$), we
have that the function 
\begin{equation*}
y(t):=\frac{d}{dt}\int_{\Omega }v(t)^{2}dx
\end{equation*}%
satisfies the inequality%
\begin{equation*}
y^{\prime }(t)+Cy^{\upsilon }(t)\leq 0
\end{equation*}%
for some $C>0$ and $\upsilon \in (0,1).$ If $\beta \in (0,1)$ then we
introduce the change of unknown $v=\mu \widehat{v}$ getting 
\begin{equation*}
\mu \widehat{v}_{t}-\mu \Delta \widehat{v}+\mu ^{\alpha }\widehat{v}^{\alpha
}=\lambda \mu ^{\beta }\widehat{v}^{\beta }.
\end{equation*}%
By choosing $\mu $ such that 
\begin{equation*}
\mu <\frac{1}{\lambda _{1}^{\frac{1}{\beta -\alpha }}}
\end{equation*}%
we can assume without loss of generality that $\lambda <\min (\lambda
_{1},1).$ Moreover, since 
\begin{equation*}
\lambda \int_{\Omega }v(t)^{\beta +1}dx\leq \lambda \int_{\Omega
}v(t)^{2}dx+\lambda \int_{\Omega }v(t)^{\alpha +1}dx,
\end{equation*}%
we get that 
\begin{equation*}
\frac{1}{2}\frac{d}{dt}\int_{\Omega }v(t)^{2}dx+\left (1-\frac{\lambda }{\lambda
_{1}}\right )\int_{\Omega }\left\vert \nabla v(t)\right\vert ^{2}dx+(1-\lambda
)\int_{\Omega }v(t)^{\alpha +1}dx\leq 0,
\end{equation*}%
and the proof ends as in the precedent case. \fin

\bigskip 

\begin{rem}
The assumption (\ref{hypo extinction}) is optimal if $\beta =1$: indeed, by
the results of \cite{D-Hern Portu} we know that for any $\lambda >\lambda
_{1}$ there exists a non-negative nontrivial solution $u$ of the associated
stationary problem $SP(\alpha ,1,\lambda )$.
\end{rem}

\bigskip 

In fact, for any $\lambda >0$ the trivial solution $u\equiv 0$ of the
stationary problem $SP(\alpha ,\beta ,\lambda )$ is asymptotically $%
\rL^{\infty }(\Omega )$-stable in the sense that it attracts solutions of $%
PP(\alpha ,\beta ,\lambda ,v_{0})$, in $\rL^{\infty }(\Omega )$,\ for small
initial data $v_{0}$.

\bigskip 

\begin{prop}
Let $v_{0}\in \rL^{\infty }(\Omega ),$ $v_{0}\geq 0.$ Assume $\beta =1$ or (%
\ref{initial datum}). Given $\lambda >0$ assume that 
\begin{equation*}
\left\Vert v_{0}\right\Vert _{{\rL}^{\infty }(\Omega )}<\lambda ^{-1/(\beta
-\alpha )}.
\end{equation*}

\noindent Then $v(t)\rightarrow 0$ in ${\rL}^{\infty }(\Omega )$ as $%
t\rightarrow +\infty .$
\end{prop}
\proof Use the solution of the associated ODE (with $\left\Vert
v_{0}\right\Vert _{{\rL}^{\infty }(\Omega )}$ as initial datum) as
supersolution. \fin

\bigskip

Concerning non-uniformly bounded trajectories we have:

\bigskip 

\begin{prop}
Let $v_{0}\in L^{\infty }(\Omega ),$ $v_{0}\geq 0$ such that 
\begin{equation}
0<u_{\lambda }(x)+\varepsilon _{0}\leq v_{0}(x)\text{ a.e. }x\in \Omega ,
\label{lower bound}
\end{equation}%
for some $\varepsilon _{0}>0$ and $u_{\lambda }$ solution of the associated
stationary problem $SP(\alpha ,\beta ,\lambda )$ such that 
\begin{equation*}
meas\{x\in \Omega :u_{\lambda }(x)+\varepsilon _{0}>\lambda ^{-1/(\beta
-\alpha )}\}>0.
\end{equation*}%
Assume $\beta =1$ or (\ref{initial datum}). Then $\left\Vert v(t)\right\Vert
_{{\rL}^{\infty }(\Omega )}\nearrow +\infty $ as $t\rightarrow +\infty .$
\end{prop}

\proof Since obviously $u_{\lambda }$ is a solution of $PP(\alpha
,\beta ,\lambda ,u_{\lambda })$ then we first get, by Theorem 2.1, that that 
$u_{\lambda }(x)\leq v(t,x)$ for any $t\in \lbrack 0,+\infty )$ and a.e. $%
x\in \Omega $. Moreover, $u_{\lambda }(x)>\lambda ^{-1/(\beta -\alpha )}>0$
on a positively measured subset $\Omega _{\lambda }$ of $\Omega $ where we
can apply the strong maximum principle to conclude that $u_{\lambda
}(x)<v(t,x)$ for any $t\in \lbrack 0,+\infty )$ and a.e. $x\in \Omega
_{\lambda }$. Since $u_{\lambda }\in \cC(\overline{\Omega })$ there exists $%
x_{\lambda }\in \overline{\Omega }_{\lambda }$ such 
\begin{equation*}
u_{\lambda }(x_{\lambda })=\min_{\overline{\Omega }_{\lambda }}u_{\lambda }
\end{equation*}%
Taking now $U(t)$ as the solution of the ODE%
\begin{equation}
ODE(\alpha ,\beta ,\lambda ,u_{\lambda }(x_{\lambda })+\varepsilon
_{0})\quad \left\{ 
\begin{array}{lr}
U_{t}+U^{\alpha }=\lambda U^{\beta } & \text{in }(0,+\infty ), \\ 
U(0)=u_{\lambda }(x_{\lambda })+\varepsilon _{0}, & 
\end{array}%
\right. 
\end{equation}%
by the standard comparison principle (notice that now the involved
nonlinearities are Lipschitz continuous on this set of values) we get that
for any $t\in \lbrack 0,+\infty )$ 
\begin{equation*}
U(t)\leq v(t,x)\text{ \ a.e. }x\in \Omega _{\lambda }.
\end{equation*}%
Finally, since we know that $U(t)\nearrow +\infty $ as $t\rightarrow +\infty 
$, we get the result. \fin

\section{Critical exponents curve on the plane ($\protect\alpha ,\protect%
\beta )$}

In this section, using Pohozaev's identity \cite{poh} and developing the
spectral analysis with respect to the fibrering procedure \cite{ilyas} we
introduce the critical exponents curve $\mathcal{C}(\rN)$ on the plane $%
(\alpha ,\beta )$ and study its main properties.

From now on we will use the notations 
\begin{eqnarray*}
T(u)=\int_\Omega |\nabla u|^2\,dx,~
A(u)=\int_\Omega |u|^{\alpha+1}\,dx,~ B(u)=\int_\Omega |u|^{\beta+1}\,dx.
\end{eqnarray*}
Then 
\begin{eqnarray}
E_\lambda(u)=\frac{1}{2}T(u) + \frac{1}{%
\alpha+1}A(u)-\lambda\frac{1}{\beta+1}B(u) .  \label{euler}
\end{eqnarray}



\noindent
\textbf{Case $0<\alpha<\beta< 1$.}

Assume that $0<\alpha <\beta <1$. Then for any fixed $u\in \rH_{0}^{1}(\Omega
)\setminus \{0\}$ the equation 
\begin{equation}
E_{\lambda }^{\prime }(ru)=0
\end{equation}%
may have at most two roots $r_{max}(u),r_{\min }(u)\in \mathbb{R}^{+}$ such
that $r_{max}(u)\leq r_{\min }(u)$. Furthermore $r_{max}(u)<r_{\min }(u)$,
if and only if 
\begin{equation*}
E_{\lambda }^{\prime \prime }(r_{max}(v)\cdot v)(v,v)<0~~\mbox{and}%
~~E_{\lambda }^{\prime \prime }(r_{\min }(v)\cdot v)(v,v)>0,
\end{equation*}%
and $r_{max}(v)=r_{\min }(v)=:r_{s}(v)$ if and only if $E_{\lambda }^{\prime
\prime }(r_{s}(v)\cdot v)=0$ (see Figure 2).

In \cite{IlEg} it has been introduced the following characteristic
(nonlinear fibrering eigenvalue) 
\begin{equation}
\Lambda _{0}=\inf_{u\in H_{0}^{1}(\Omega )\setminus 0}\lambda _{0}(u),
\end{equation}%
where 
\begin{equation*}
\lambda _{0}(u)=c_{0}^{\alpha ,\beta }\lambda (u),
\end{equation*}%
\begin{equation*}
c_{0}^{\alpha ,\beta }=\frac{(1-\alpha )(1+\beta )}{(1-\beta )(1+\alpha )}%
\left( \frac{(1+\alpha )(1-\beta )}{2(\beta -\alpha )}\right) ^{\frac{\beta
-\alpha }{1-\alpha }}
\end{equation*}%
and 
\begin{equation}
\lambda (u)=\frac{A(u)^{\frac{1-\beta }{1-\alpha }}T(u)^{\frac{\beta -\alpha 
}{1-\alpha }}}{B(u)}.  \label{lam}
\end{equation}%
Note that by the Gagliardo-Nirenberg inequality (see \cite[Proposition 2]%
{IlEg}) it follows that $0<\Lambda _{0}<+\infty $. In~\cite{IlEg}, it was
proved the

\begin{prop}
\label{lam0} If $\lambda \geq \Lambda _{0}$, then there exists $u\in
H_{0}^{1}(\Omega )\setminus \{0\}$ such that $E_{\lambda }^{\prime }(u)=0$ and $%
E_{\lambda }(u)\leq 0$, $E_{\lambda }^{\prime \prime }(u)>0$.
\end{prop}

We need also the following characteristic value from \cite{IlEg} 
\begin{equation}
\Lambda _{1}=\inf_{u\in \rH_{0}^{1}\setminus \{0\}}\lambda _{1}(u).  \label{P30}
\end{equation}%
where 
\begin{equation}
\lambda _{1}(u)=c_{1}^{\alpha ,\beta }\lambda (u),  \label{P20}
\end{equation}%
where 
\begin{equation}
c_{1}^{\alpha ,\beta }=\frac{1-\alpha }{1-\beta }\left( \frac{1-\beta }{%
\beta -\alpha }\right) ^{\frac{\beta -\alpha }{1-\alpha }}.  \label{P40}
\end{equation}%
As before we have $0<\Lambda _{1}<+\infty $. Furthermore, $0<\Lambda
_{1}<\Lambda _{0}<+\infty $ (see \cite[Claim 2]{IlEg}) and we have as in
Lemma \ref{lam0} (see also \cite{IlEg})

\begin{prop}
\label{lam1} If $\lambda >\Lambda _{1}$, then there exists $u\in
\rH_{0}^{1}(\Omega )\setminus \{0\}$ such that $E_{\lambda }^{\prime }(u)=0$,
whereas if $\lambda <\Lambda _{1}$, then $E_{\lambda }^{\prime }(u)>0$ for
any $u\in \rH_{0}^{1}(\Omega )\setminus \{0\}$.
\end{prop}

Let $u\in\rH_{0}^{1}(\Omega )$ be a weak solution of \eqref{1}. Standard
regularity arguments show that $u\in \cC^{1,\gamma }(\overline{\Omega })\cap
\cC^{2}(\Omega )$ for some $\gamma \in (0,1)$. Note that by the assumption $%
\partial \Omega $ is a $\cC^{1}$-manifold. Therefore Pohozaev's identity holds 
\cite{poh,Takac_Ilyasov}, namely 
\begin{equation}
P_{\lambda }(u)+\frac{1}{2\rN}\int_{\partial \Omega }\left\vert \frac{\partial
u}{\partial \nu }\right\vert ^{2}\,x\cdot \nu \,ds=0,  \label{poh}
\end{equation}%
where 
\begin{equation*}
P_{\lambda }(u):=\displaystyle{\frac{\rN-2}{2\rN}T(u)+\frac{1}{\alpha +1}A(u)}-\lambda \frac{1}{\beta +1%
}B(u),~~~u\in \rH_{0}^{1}(\Omega ).
\end{equation*}

Note that if $\Omega $ is a star-shaped (strictly star-shaped) domain with
respect to the origin of $\mathbb{R}^{\rN}$, then $x\cdot \nu \geq 0$ ($x\cdot
\nu >0$) for all $x\in \partial \Omega $. Thus we have

\begin{prop}
\label{PropP} Assume that $\Omega $ is a star-shaped domain with respect to
the origin of $\mathbb{R}^{\rN}$, then $P_{\lambda }(u)\leq 0$ ($P_{\lambda
}(u)=0$) for any weak (flat or compactly supported) solution $u$ of \eqref{1}%
. If, in addition, $\Omega $ is strictly star-shaped, then a weak solution $%
u $ of \eqref{1} is flat or it has compact support if and only if $%
P_{\lambda }(u)=0$.
\end{prop}

Let us study the critical exponent curve $\mathcal{C}(\rN)$ (see \eqref{crN})
and prove Lemma \ref{lemG}. Consider the system (see~\cite{IlEg}) 
\begin{equation}
\left\{ 
\begin{array}{l}
E_{\lambda }^{\prime }(u):=\ T(u)+A(u)-\lambda B(u)=0 \\ 
\\ 
P_{\lambda }(u):=\displaystyle{\frac{\rN-2}{2\rN}T(u)+\frac{1}{\alpha +1}A(u)-\lambda \frac{1}{\beta +1
}B(u)=0} \\ 
\\ 
E_{\lambda }^{\prime \prime }(u):=T(u)+\alpha A(u)-\lambda \beta B(u)=0.%
\end{array}%
\right.  \label{sis1}
\end{equation}%
This system is solvable with respect to the variables $T(u),A(u),B(u)$ if
the corresponding determinant 
\begin{equation}
D=\frac{(\beta -\alpha )(2(1+\alpha )(1+\beta )-\rN(1-\alpha )(1-\beta ))}{%
2\rN(1+\alpha )(1+\beta )}.  \label{theta}
\end{equation}%
is non-zero.

\noindent On the other hand $D=0$ if and only if $(\alpha ,\beta )\in 
\mathcal{C}(\rN)$.

\bigskip

\noindent {\sc Proof of Lemma \ref{lemG}}. \quad Let $\Omega $ be a star-shaped
domain with respect to the origin of $\mathbb{R}^{\rN}$. Then by Proposition %
\ref{PropP} we have $P_{\lambda }(u)=0$ for any flat or compactly supported
solution $u$ of \eqref{1}. Note also that $E_{\lambda }^{\prime }(u)=0$.
Thus, in case $(\alpha ,\beta )\in \mathcal{C}(\rN)$, i.e. when the
determinant of system \eqref{sis1} is equal to zero one has $E_{\lambda
}^{\prime \prime }(u)=0$ and we get the proof of statement 1), Lemma \ref%
{lemG}. Observe 
$$
\begin{array}{ll}
D\cdot \frac{2\rN(1+\alpha )}{(1-\alpha )[-2(1+\alpha )-\rN(1-\alpha )]} B(u)=&
\displaystyle \frac{1}{1-\alpha }(E_{\lambda }^{\prime \prime }(u)-E_{\lambda }^{\prime
}(u))-  \label{D1} \\ [.2cm]
& 
\displaystyle \frac{2N(1+\alpha )}{(\rN-2)(1+\alpha )-2\rN}(P_{\lambda }(u)-\frac{N-2}{2N}%
E_{\lambda }^{\prime }(u)).  \notag
\end{array}%
$$
Thus if $(\alpha ,\beta )\in \mathcal{E}_{u}(\rN)$ and $P_{\lambda }(u)=0$, $%
E_{\lambda }^{\prime }(u)=0$, then 
\begin{equation*}
E_{\lambda }^{\prime \prime }(u)=-D\cdot \frac{2\rN(1+\alpha )}{(1-\alpha
)[2(1+\alpha )+\rN(1-\alpha )]}B(u)<0
\end{equation*}%
and we obtain the proof of statement 2), Lemma \ref{lemG}.

Under assumption 3) of Lemma \ref{lemG}, for a weak solution $u$ of \eqref{1}
we have $P_{\lambda }(u)\leq 0$ (see Proposition \ref{PropP}) and therefore %
\eqref{D1} yields 
\begin{equation*}
E_{\lambda }^{\prime \prime }(u)\geq -D\cdot \frac{2\rN(1+\alpha )}{(1-\alpha
)[-2(1+\alpha )-\rN(1-\alpha )]}B(u)>0,
\end{equation*}%
since $D>0$ for $(\alpha ,\beta )\in \mathcal{E}_{s}(\rN)$. This completes the
proof of Lemma \ref{lemG}. \fin

\bigskip

\noindent
\textbf{Case $\beta=1$.}

Recall some results from \cite{diaz-H-Il}. In what follows $(\lambda
_{1},\varphi _{1})$ denotes the first eigenpair of the operator $-\Delta $
in $\Omega $ with zero boundary conditions. Let $u\in \rH_{0}^{1}(\Omega )$. The fibrering mapping in this case is
defined by 
\begin{equation*}
\Phi _{u}(r)=E_{\lambda }(ru)=\frac{r^{2}}{2}\rH_{\lambda }(u)+\frac{%
r^{1+\alpha }}{1+\alpha }A(u)
\end{equation*}%
where we denote 
\begin{equation*}
\rH_{\lambda }(u):=\int_{\Omega }|\nabla u|^{2}\,\mathrm{d}x-\lambda
\int_{\Omega }|u|^{2}\,\mathrm{d}x.
\end{equation*}%
Then 
\begin{equation*}
\Phi _{u}^{\prime }(r)=E_{\lambda }^{\prime }(ru)=r\rH_{\lambda }(u)+r^{\alpha
}A(u)
\end{equation*}%
and the equation $\Phi _{u}^{\prime }(r)=0$ has a positive solution only if
both term in $\Phi _{u}^{\prime }(r)$ have opposite sign, that is if and
only if $\rH_{\lambda }(u)<0$. Note that there is $u\in \rH_{0}^{1}(\Omega )$
such that $\rH_{\lambda }(u)<0$ iff $\lambda >\lambda _{1}$. It turns out that
the only point $r(u)$ where $\Phi _{u}^{\prime }(r)=0$ is given by 
\begin{equation}
r(u)=\left( \frac{A(u)}{-\rH_{\lambda }(u)}\right) ^{1/(1-\alpha )}.  \label{t}
\end{equation}%
Furthermore, $E_{\lambda }^{\prime \prime }(r(u)u)(u,u)<0$ and 
\begin{equation}
E_{\lambda }(r(u)u)=\max_{r>0}E_{\lambda }(ru)  \label{maxr}.
\end{equation}%
Substituting \eqref{t} into $E_{\lambda }(ru)$ we obtain 
\begin{equation}
J_{\lambda }(u):=E_{\lambda }(r_{\lambda }(u)u)=\frac{(1-\alpha )}{%
2(1+\alpha )}\frac{A(u)^{\frac{2}{1-\alpha }}}{(-\rH_{\lambda }(u))^{\frac{%
1+\alpha }{1-\alpha }}}.  \label{JL}
\end{equation}%
Consider 
\begin{equation}
\widehat{E}_{\lambda }=\min \{J_{\lambda }(u):~~u\in \rH_{0}^{1}(\Omega )\setminus
\{0\},~\rH_{\lambda }(u)<0\}.  \label{minJomega}
\end{equation}%
It follows directly

\begin{prop}
\label{prop:fc} A point $u\in \rH_{0}^{1}(\Omega )$ is a minimizer of %
\eqref{minJomega} if and only if $\tilde{u}=r(u)u$ is a ground state of~\eqref{1R}.
\end{prop}

\bigskip

\begin{rem}
We point out that in both cases, $\beta <1$ and $\beta =1,$ the above
results can be extended to the case in which the ground solution of $%
SP(\alpha ,\beta ,\lambda )$ minimizes the energy on the closed convex cone 
\begin{equation*}
K=\{v\in \rH_{0}^{1}(\Omega ),v\geq 0\text{ on }\Omega \}.
\end{equation*}%
Indeed, we introduce the modified energy functional 
\begin{equation*}
E_{\lambda }^{+}(u)=E_{\lambda }(u)+\int_{\Omega }j(u)dx
\end{equation*}%
where%
\begin{equation*}
j(u)=\left\{ 
\begin{array}{cr}
0 & \text{if }u\in K, \\ 
+\infty & \text{otherwise.}%
\end{array}%
\right.
\end{equation*}%
Notice that $j(ru)=j(u)$ for any $r>0.$ Obviously $E_{\lambda
}^{+}(u)=E_{\lambda }(u)$ if $u\in K.$ Moreover the additional term arising
in the associated Euler-Lagrange equation, given by the subdifferential of
the convex function $\displaystyle \int_{\Omega }j(u)dx$, vanishes when the ground state
solution of $SP(\alpha ,\beta ,\lambda )$ is nonnegative.
\end{rem}

\bigskip

\section{Existence of ground state}

In this Section, we prove the first parts of Theorems \ref{Th1}, \ref{Th1.1}.

{\sc Proof of (1), Theorem \ref{Th1}}\quad Assume $\beta <1$. In this case, the existence of a ground state
of \eqref{1} when $(\alpha ,\beta )\in \mathcal{E}_{s}(\rN)$ has been proved
in \cite{IlEg}. The proof for the points $(\alpha ,\beta )\in \mathcal{E}%
\setminus \mathcal{E}_{s}(\rN)$ can be obtained in a similar way. However for
the sake of completeness, we present a summary of the proof.

Consider the constrained minimization problem of $E_{\lambda }(u)$ on the
associated Nehari manifold 
\begin{equation}
\left\{ 
\begin{array}{l}
\ E_{\lambda }(u)\rightarrow \min \\ 
\\ 
E_{\lambda }^{\prime }(u)(u)=0.%
\end{array}%
\right.  \label{min1}
\end{equation}%
We denote by 
\begin{equation*}
\mathcal{N}_{\lambda }:=\{u\in H_{0}^{1}(\Omega ):~E_{\lambda }^{\prime
}(u)=0\}
\end{equation*}%
the admissible set of (\ref{min1}), i.e. the corresponding Nehari manifold.
Denote also 
\begin{equation*}
\widehat{E_{\lambda }}:=\min \{E_{\lambda }(u):~u\in \mathcal{N}_{\lambda }\}
\end{equation*}%
the minimum value in this problem. Note that by Proposition \ref{lam1}, 
$\cN_{\lambda }\neq \emptyset $ for any $\lambda >\Lambda _{1}$.
Furthermore, by Sobolev's inequalities we have 
\begin{equation*}
E_{\lambda }(u)\geq \frac{1}{2}||u||_{1}^{2}-c_{1}||u||_{1}^{1+\beta
}\rightarrow \infty
\end{equation*}%
as $||u||_{1}\rightarrow \infty $, since $2>1+\beta $. Thus $E_{\lambda }(u)$
is a coercive functional on $\rH_{0}^{1}(\Omega )$. Using this it is not hard
to prove the following (see also \cite[Lemma 9]{IlEg})

\begin{prop}
\label{le1e} Let $(\alpha ,\beta )\in \mathcal{E}$. Then for any $\lambda
\geq \Lambda _{1}$ problem (\ref{min1}) has a minimizer $u_{\lambda }\in
\rH_{0}^{1}(\Omega )\setminus \{0\}$, i.e. $E_{\lambda }(u_{\lambda })=\widehat{%
E_{\lambda }}$ and $u_{\lambda }\in \cN_{\lambda }$.
\end{prop}

Let $\lambda \geq \Lambda _{1}$ and $u_{\lambda }\in \rH_{0}^{1}(\Omega
)\setminus \{0\}$ be a minimizer of (\ref{min1}). Then by the Lagrange
multipliers rule there exist $\mu _{1}$, $\mu _{2}$ such that 
\begin{equation}
\mu _{1}DE_{\lambda }(u_{\lambda })=\mu _{2}DE_{\lambda }^{\prime
}(u_{\lambda })(u_{\lambda }),  \label{eq2}
\end{equation}%
and $|\mu _{1}|+|\mu _{2}|\neq 0$. Thus, if $\mu _{2}=0$, then $u_{\lambda }$
is a weak solution of (\ref{1}).

This condition is satisfied under the assumptions of the following result:

\begin{prop}
\label{Lag2} Let $(\alpha ,\beta )\in \mathcal{E}$. Then for any $\lambda
\geq \Lambda _{0}$ (\ref{1}) has a ground state $u_{\lambda }$ which is
nonnegative, $u\in \cC^{1,\gamma }(\overline{\Omega })\cap \cC^{2}(\Omega )$ for
some $\gamma \in (0,1)$ and $E_{\lambda }^{\prime \prime }(u_{\lambda
})(u_{\lambda },u_{\lambda })>0$.
\end{prop}

\proof Since $0<\Lambda _{1}<\Lambda _{0}$, then by
Proposition \ref{le1e} for any$\lambda \geq \Lambda _{0}$ there exists a
minimizer $u_{\lambda }\in \rH_{0}^{1}(\Omega )\setminus \{0\}$ of (\ref{min1}).
Lemma \ref{lam0} implies that there is $u\in \mathcal{N}_{\lambda }$ such
that $E_{\lambda }(u)\leq 0$ and therefore $E_{\lambda }(u_{\lambda })\leq
E_{\lambda }(u)\leq 0$. This implies that $E_{\lambda }^{\prime \prime
}(u_{\lambda })(u_{\lambda },u_{\lambda })>0$. Let us test \eqref{eq2} by $%
u_{\lambda }$. Then 
\begin{equation*}
\mu _{1}E_{\lambda }^{\prime }(u_{\lambda })(u_{\lambda })=\mu
_{2}(E_{\lambda }^{\prime \prime }(u_{\lambda })(u_{\lambda },u_{\lambda
})+E_{\lambda }^{\prime }(u_{\lambda })(u_{\lambda })).
\end{equation*}%
Since $E_{\lambda }^{\prime }(u_{\lambda })(u_{\lambda })=0$, this yields
that $\mu _{2}E_{\lambda }^{\prime \prime }(u_{\lambda })=0$. But $%
E_{\lambda }^{\prime \prime }(u_{\lambda })(u_{\lambda },u_{\lambda })\neq 0$
and therefore $\mu _{2}=0$. Thus, by \eqref{eq2} we obtain $DE_{\lambda
}(u_{\lambda })=0$, i.e $u_{\lambda }$ is a weak solution of (\ref{1}).
Since any weak solution $w_{\lambda }$ of (\ref{1}) belongs to $\cN%
_{\lambda }$, then \eqref{min1} yields that $u_{\lambda }$ is a ground
state. The rest of the lemma is proved by standard way.~\fin

\bigskip

From this Proposition arguing by contradiction, it is not hard to show that
there is an interval $(\Lambda _{0}-\varepsilon ,+\infty )$ for some $%
\varepsilon >0$ such that for any $\lambda \in (\Lambda _{0}-\varepsilon
,+\infty )$ the minimizer $u_{\lambda }$ of (\ref{min1}) satisfies $%
E_{\lambda }^{\prime \prime }(u_{\lambda })>0$. From this, as in the proof
of Proposition \ref{Lag2}, it follows that $u_{\lambda }$ is a ground state
of (\ref{1}) which is nonnegative and $u\in \cC^{1,\gamma }(\overline{\Omega }%
)\cap \cC^{2}(\Omega )$ for some $\gamma \in (0,1)$.

Thus we have a proof that there exists $\lambda ^{\ast }\in (\Lambda
_{1},\Lambda _{0})$ such that for all $\lambda >\lambda ^{\ast }$ problem (%
\ref{1}) has a ground state $u_{\lambda }$, which is nonnegative in $\Omega $%
, $u\in \cC^{1,\gamma }(\overline{\Omega })\cap \cC^{2}(\Omega )$ for some $%
\gamma \in (0,1)$ and $E_{\lambda }^{\prime \prime }(u_{\lambda
})(u_{\lambda },u_{\lambda })>0$ This completes the proof of statement (1)
of Theorem \ref{Th1}.
\par
\noindent 
{\sc Proof of (1), Theorem \ref{Th1.1}}\quad The existence of a ground state is obtained from the constrained
minimization problem \eqref{minJomega} and then using Proposition \ref%
{prop:fc}. The implementation of this proof was done in \cite[Theorem 2.1,
p.6]{diaz-H-Il}.

%
%
%

\section{Existence of ground state flat solutions in case $\protect\beta=1$}

In this Section, we prove statement (1) in Theorem \ref{Th2.2}. Consider now the following auxiliary problem on the whole space $\mathbb{R}%
^{\rN}$: 
\begin{equation}
\left\{ 
\begin{array}{ll}
-\Delta u+u^{\alpha }=u & \text{in }\mathbb{R}^{\rN}\text{,} \\ 
\ u\geq 0 & \text{on }\mathbb{R}^{\rN}.%
\end{array}%
\right.  \label{1R}
\end{equation}%
Here and subsequently, $\rH^{1}(\mathbb{R}^{\rN})$ denotes the standard Sobolev
space with the norm 
\begin{equation*}
||u||_{1}=\left (\int_{\mathbb{R}^{\rN}}|u|^{2}\,dx+\int_{\mathbb{R}^{\rN}}|\nabla
u|^{2}\,dx\right )^{1/2}.
\end{equation*}

Problem (\ref{1R}) has a variational form with the Euler-Lagrange functional 
\begin{equation*}
E(u)=\frac{1}{2}\rH(u)+\frac{1}{\alpha +1}A(u),\quad u\in \rW^{1,2}(\mathbb{R}%
^{\rN})
\end{equation*}%
where 
\begin{equation*}
\rH(u)=\int_{\mathbb{R}^{\rN}}|\nabla u|^{2}\,dx-\int_{\mathbb{R}%
^{\rN}}|u|^{2}\,dx,~A(u)=\int_{\mathbb{R}^{\rN}}|u|^{\alpha +1}\,dx.
\end{equation*}%
As above we call a nonzero weak solution $u_{\lambda }$ of \eqref{1R} a
ground state of \eqref{1R} if it holds 
\begin{equation*}
E(u_{\lambda })\leq E(w_{\lambda })
\end{equation*}%
for any nonzero weak solution $w_{\lambda }$ of \eqref{1R}. 
The fibreing map in this case is given as follows 
\begin{equation*}
\Phi _{u}(r):=E(ru)=\frac{r^{2}}{2}\rH(u)+\frac{r^{1+\alpha }}{\alpha +1}%
A(u),~~~u\in \rH^{1}(\mathbb{R}^{\rN}),~~t\in \mathbb{R}^{+}
\end{equation*}%
and, for fix $u\in \rH^{1}(\mathbb{R}^{\rN})$ the equation 
\begin{equation*}
\Phi _{u}^{\prime }(r)\equiv r\rH(u)+r^{\alpha }A(u)=0,\quad r\in \mathbb{R}%
^{+}.
\end{equation*}%
has only one root 
\begin{equation}
r(u)=\left( \frac{A(u)}{-\rH_{\lambda }(u)}\right) ^{1/(1-\alpha )}  \label{th}
\end{equation}%
which exists if and only if $\rH(u)<0$.

As above, substituting this root into $E_{\lambda }(ru)$ we obtain a
zero-homogeneous functional 
\begin{equation}
J(u):=E(r(u)u)=\frac{(1-\alpha )}{2(1+\alpha )}\frac{A(u)^{\frac{2}{1-\alpha 
}}}{(-\rH(u))^{\frac{1+\alpha }{1-\alpha }}},  \label{fibF}
\end{equation}%
and we consider 
\begin{equation}
\widehat{E}^{\infty }=\min \{J(u):~~u\in \rH^{1}(\mathbb{R}^{\rN})\setminus
\{0\},~\rH(u)<0\}.  \label{minJinfty}
\end{equation}%
As above, it follows directly the

\begin{prop}
\label{prop:fcinfty}We have that $u$ is a minimizer of \eqref{minJinfty} if
and only if $\tilde{u}=r(u)u$ is a ground state of \eqref{1R}.
\end{prop}

In Appendix below, using \eqref{minJinfty} we prove the

\begin{lemma}
\label{lem:r1} Assume $0<\alpha <1$. Then problem \eqref{1R} has a classical
nonnegative solution $u\in \rH^{1}(\mathbb{R}^{\rN})$ which is a ground state.
\end{lemma}

The following result can be found in \cite{Serrin-Zou}

\begin{lemma}
\label{lem:3} Assume $0<\alpha <1$. Then any classical solution $u$ of %
\eqref{1R} has a compact support. Furthermore if we define 
\begin{equation*}
\Theta :=\{x\in \mathbb{R}^{\rN}:u(x)>0\}.
\end{equation*}%
Then for every connected component $\Xi $ of $\Theta $ we have
\begin{enumerate}
\item $\Xi$ is a ball;
\item $u$ is radially symmetric with respect to the centre of the ball $\Xi $.
\end{enumerate}
\end{lemma}

Lemmas \ref{lem:r1}, \ref{lem:3} yield

\begin{coro}
\label{corTR} Assume $0<\alpha<1$. Then there is a radius $R^*>0$ such that
problem \eqref{1R} has a ground state $u^*$ which is a flat classical radial
solution and 
\begin{equation*}
\mathrm{supp}(u^*)=B_{R^*}.
\end{equation*}
\end{coro}

Let us return to problem \eqref{1}. From Corollary \ref{corTR} we have

\begin{coro}
\label{corTR2} Assume that $B_{R^*} \subset \Omega$. Then the ground state $%
u_\lambda$ of \eqref{1} with $\lambda =1$ coincides with the ground state $%
u^*$ of \eqref{1R} that is $u_\lambda|_{\lambda=1}$ is a compact support
classical radial solution and 
\begin{equation*}
\mathrm{supp}(u_\lambda)|_{\lambda=1}\equiv \bar{\Theta}=B_{R^*}.
\end{equation*}
\end{coro}

\proof Any function $w$ from $\rH_{0}^{1}(\Omega )$
can be extended to $\mathbb{R}^{\rN}$ as 
\begin{equation}
\left\{ 
\begin{array}{ll}
\tilde{w}=w &\mbox{in}~\Omega , \\ 
\tilde{w}=0 & \mbox{in}~\mathbb{R}^{\rN}\setminus \Omega ,%
\end{array}%
\right.  \label{expan}
\end{equation}%
Then $\tilde{w}\in \rH^{1}(\mathbb{R}^{\rN})$ and in this sense we may assume
that $\rH_{0}^{1}(\Omega )\subset \rH^{1}(\mathbb{R}^{\rN})$. Therefore 
\begin{equation*}
\widehat{E}^{\infty }\leq \widehat{E}_{1}\equiv \min \{J_{1}(v):~~v\in
\rH_{0}^{1}(\Omega )\setminus 0,v\geq 0,~\rH_{1}(v)<0\}.
\end{equation*}%
Note that $u^{\ast }\in K\subset \rH_{0}^{1}(B_{R^{\ast }})\subset
\rH_{0}^{1}(\Omega )$. This yields $\widehat{E}^{\infty }=E(u^{\ast })=\widehat{E}_{1}$
and we get the proof. \fin

\bigskip

Assume now that $\Omega $ is a is star-shaped domain in $\mathbb{R}^{\rN}$,
with respect to the some point $z\in \mathbb{R}^{\rN}$ which without loss of
generality we may assume coincides with the origin $0\in \mathbb{R}^{\rN}$.

Let $u_{\lambda }$ be a ground state of \eqref{1}. By making a change of
variable $v_{\lambda (\kappa )}(y)=\kappa ^{-2/(1-\alpha )}u_{\lambda
}(\kappa y),~~y\in \Omega _{\kappa }$, with $\kappa >0$ we get 
\begin{equation}
\left\{ 
\begin{array}{ll}
\ -\Delta v_{\lambda (\kappa )}=\lambda (\kappa )v_{\lambda (\kappa
)}-v_{\lambda (\kappa )}^{\alpha } &\mbox{in}~\Omega _{\kappa }, \\ [.175cm]
\ v_{\lambda (\kappa )}=0 & \mbox{on}~\Omega _{\kappa }.%
\end{array}%
\right.  \label{1r2}
\end{equation}%
where $\lambda (\kappa )=\lambda \kappa ^{2}$, $\Omega _{\kappa }=\{y\in 
\mathbb{R}^{N}:~y=x/\kappa ,~x\in \Omega \}$. Since $u_{\lambda }$ is a
ground state of \eqref{1}, then it is easy to see that $v_{\lambda (\kappa
)} $ is also a ground state of \eqref{1r2}. Note that if $\kappa =\sqrt{%
1/\lambda }$ then $\lambda (\kappa )=1$. On the other hand, if $\kappa $ is
sufficiently small then $B_{R^{\ast }}\subset \Omega _{\kappa }$. Hence by
Corollary \ref{corTR} there is a sufficiently large $\lambda ^{\ast }$ such
that for any $\lambda >\lambda ^{\ast }$ the ground state $v_{\lambda
(\kappa )}$ with $\lambda (\kappa )=\lambda \cdot (\kappa )^{2}$, $\kappa =%
\sqrt{1/\lambda }$ is a flat or compactly supported classical radial
solution of \eqref{1r2} which coincides with the ground state $u^{\ast }$ of %
\eqref{1R}. Thus we have proved

\begin{coro}
\label{corGB} Assume $0<\alpha<1$. Then there exists $\lambda^*>0$ such that
for any $\lambda\geq \lambda^*$ problem \eqref{1} has a ground state $%
u_{\lambda}$ which is a flat classical radial solution. Furthermore, $%
u_{\lambda^*}(x)= \kappa^{2/(1-\alpha)}u^* (x/\kappa)$ where $\kappa=\sqrt{%
1/\lambda}$ and $u^*$ is a flat classical radial ground state of \eqref{1R}.
\end{coro}

Note that by \cite[Lemma 3.3]{diaz-H-Il} 
\begin{equation*}
\lambda ^{\ast }>\lambda ^{c}=(1+\frac{2(1+\alpha )}{\rN(1-\alpha )})\cdot
\lambda _{1}(\Omega ).
\end{equation*}%
Furthermore, for any $\lambda \in (\lambda _{1}(\Omega ),\lambda ^{c})$
problem \eqref{1} cannot have flat solutions in $C^{1}(\overline{\Omega })$.

\section{Lyapunov stability of flat ground states}

In this Section, first we prove statements (2) of Theorem \ref{Th1} and then
prove \textbf{(III)} of Theorem \ref{Th2} .

To prove the stability we will use the Lyapunov Function method. Let $%
u_{\lambda }$ be a ground state of \eqref{1} such that $E_{\lambda }^{\prime
\prime }(u_{\lambda })(u_{\lambda },u_{\lambda })>0$. For $\delta >0$,
denote 
\begin{equation*}
U_{\delta }(u_{\lambda }):=\{v\in \rH_{0}^{1}(\Omega ):~||u_{\lambda
}-v||<\delta \}.
\end{equation*}%
Observe that $E_{\lambda },E_{\lambda }^{\prime \prime }:\rH_{0}^{1}(\Omega
)\rightarrow \mathbb{R}$ are continuous maps. Hence there exists $\delta
_{0}>0$ such that $E_{\lambda }^{\prime \prime }(u)(u,u)>0$ for all $u\in
U_{\delta }(u_{\lambda })$ if $0<\delta <\delta _{0}$.

In the next two lemmas we show that $E_{\lambda }$ is a Lyapunov function in
the neighborhood $U_{\delta }(u_{\lambda })$ if $0<\delta <\delta _{0}$.

\begin{lemma}
\label{lemSt1} Assume {\rm ({\bf U})}. Let $\lambda >\lambda ^{\ast }$ and $u_{\lambda }$
be a ground state of \eqref{1} such that $E_{\lambda }^{\prime \prime
}(u_{\lambda })>0$. Then for any $\delta \in (0,\delta _{0})$ it satisfies 
\begin{equation}
E_{\lambda }(u)>E_{\lambda }(u_{\lambda })=\hat{E_{\lambda }}~~\forall u\in
U_{\delta }(u_{\lambda })\setminus \{u_{\lambda }\}
\end{equation}
\end{lemma}

\proof Suppose, contrary to our claim that for every $%
\delta \in (0,\delta _{0})$ there exists $u^{\delta }\in U_{\delta
}(u_{\lambda })\setminus \{u_{\lambda}\} $ such that $E_{\lambda }(u^{\delta
})\leq E_{\lambda }(u_{\lambda })$. This implies that there exists a
sequence $u^{n}\in U_{\delta _{0}}(u_{\lambda })$ such that $%
u^{n}\rightarrow u_{\lambda }$ in $\rH_{0}^{1}(\Omega )$ as $n\rightarrow
\infty $ and 
\begin{equation}
E_{\lambda }(u^{n})\leq E_{\lambda }(u_{\lambda })\quad n=1,2,....  \label{st1}
\end{equation}%
Note that by property \textbf{(U)} we may assume that the point $u^{n}$ for
any $n=1,2,...,$ is not a ground state of \eqref{1}. Furthermore, $r_{\min
}(u_{\lambda })=1$ since $E_{\lambda }^{\prime \prime }(u_{\lambda })>0$.
Thus by \eqref{min1} we have 
\begin{equation*}
E_{\lambda }(r_{miN}(u^{n})u^{n})>E_{\lambda }(u_{\lambda })\quad n=1,2,....
\end{equation*}%
Moreover, this and \eqref{st1} yield that 
\begin{equation}
1<r_{max}(u^{n})<r_{\min }(u^{n}).  \label{st3}
\end{equation}%
Note that $r_{max}(\cdot ),r_{\min }(\cdot ):\rH_{0}^{1}(\Omega )\rightarrow 
\mathbb{R}$ are continuous maps. Hence 
\begin{equation*}
r_{\min }(u^{n})\rightarrow r_{\min }(u_{\lambda })=1~~\mbox{as }%
n\rightarrow \infty ,
\end{equation*}%
since $u^{n}\rightarrow u_{\lambda }$ in $\rH_{0}^{1}(\Omega )$ as $%
n\rightarrow \infty $. Then by \eqref{st3} we have also 
\begin{equation*}
r_{max}(u^{n})\rightarrow r_{\min }(u_{\lambda })=1~~\mbox{as}\quad n\rightarrow
\infty .
\end{equation*}%
From this and since $E_{\lambda }^{\prime \prime }(r_{max}(u^{n})u^{n})\leq
0 $ and $E_{\lambda }^{\prime \prime }(r_{\min }(u^{n})u^{n})\geq 0$ we
conclude that 
\begin{equation*}
E_{\lambda }^{\prime \prime }(u_{\lambda })=0.
\end{equation*}%
But this is impossible by the assumption. This contradiction completes the
proof. \fin

\begin{lemma}
Let $v(t)$, $t \in [0,T)$ be a weak solution of \eqref{p1}. Then 
\begin{equation}
\frac{\partial}{\partial t} E_\lambda(v(t)) \leq 0 ~~~\mbox{in} ~~~ (0,T).
\end{equation}
\end{lemma}

\textit{Proof.}\thinspace\ By the additional regularity obtained in Section
2, there exists $\frac{\partial }{\partial t}E_{\lambda }(v(t))$ in $(0,T)$
and 
$$
\frac{\partial }{\partial t}E_{\lambda }(v(t))=D_{u}E_{\lambda
}(v(t))(v_{t}(t))= 
<-\Delta v(t)-\lambda |v|^{\beta -1}v+|v|^{\alpha
-1}v,v_{t}(t)>=-||v_{t}(t)||_{L^{2}}^{2}\leq 0.
$$
Thus we get the result.\fin

\bigskip

The proof of (2), Theorem \ref{Th1} will follow from

\begin{lemma}
\label{lemst2} Assume (U). Let $\lambda >\lambda ^{\ast }$ and $u_{\lambda }$
be a ground state of \eqref{1} such that $E_{\lambda }^{\prime \prime
}(u_{\lambda })>0$. Then for any given $\varepsilon >0$, there exists $%
\delta \in (0,\delta _{0})$ such that 
\begin{equation}
||u_{\lambda }-v(t;w_{0})||_{1}<\varepsilon ~~\mbox{for any}~w_{0}\geq 0~%
\text{such that }~~||u_{\lambda }-w_{0}||_{1}<\delta ,~~\forall t>0.
\end{equation}
\end{lemma}

\proof Without loss of generality we may assume that $%
\varepsilon \in (0,\delta _{0})$. Consider 
\begin{equation}
d_{\varepsilon }:=\inf \{E_{\lambda }(w):w\in \rH_{0}^{1}(\Omega
),||u_{\lambda }-w||_{1}=\varepsilon \}.
\end{equation}%
Then $d_{\varepsilon }>\hat{E_{\lambda }}$. Indeed, assume the opposite,
that there is a sequence $w^{n}\in K$, $||u_{\lambda
}-w^{n}||_{1}=\varepsilon $ and $E_{\lambda }(w^{n})\rightarrow \hat{%
E_{\lambda }}$. Hence $(w^{n})$ is bounded in $\rH_{0}^{1}(\Omega )$ and
therefore by the embedding theorem there exists a subsequence (again denoted
by $(w^{n})$) such that $w^{n}\rightarrow w_{0}$ weakly in $\rH_{0}^{1}(\Omega
)$ and strongly in $\rL_{p},~ 1<p<2^{\ast }$ for some $w_{0}\in
\rH_{0}^{1}(\Omega )$. Since $||u||_{1}^{2}$ is a weakly lower semi-continuous
functional on $\rH_{0}^{1}(\Omega )$, one has $\hat{E_{\lambda }}\geq
E_{\lambda }(w_{0})$ and $||u_{\lambda }-w_{0}||_{1}\leq \varepsilon $. By
Lemma \ref{lemSt1} this is possible only if $w_{0}$ is a ground state of %
\eqref{1}, i.e., a minimizer of \eqref{min1}. But then $\widehat{E_{\lambda }}%
=E_{\lambda }(w_{0})$ implies that $w^{n}\rightarrow w_{0}$ strongly in $%
\rH_{0}^{1}(\Omega )$. From here we have $\varepsilon =||u_{\lambda
}-w^{n}||_{1}\rightarrow ||u_{\lambda }-w_{0}||_{1}$. Thus $w_{0}\in
U_{\delta _{0}}(u_{\lambda })$ and $u_{\lambda }\neq w_{0}$. Since by
property \textbf{(U)} $u_{\lambda }$ is the unique non-negative solution of %
\eqref{1} in $U_{\delta _{0}}(u_{\lambda })$ we get a contradiction.

Let $\sigma >0$ be an arbitrary value such that $d_{\varepsilon }-\sigma >%
\widehat{E_{\lambda }}$. Then by continuity of $E_{\lambda }(w)$ one can find $%
\delta \in (0,\varepsilon )$ such that 
\begin{equation}
E_{\lambda }(w)<d_{\varepsilon }-\sigma \quad\forall w\in U_{\delta
}(u_{\lambda })\subset U_{\varepsilon }(u_{\lambda }).  \label{eqst3}
\end{equation}%
We claim that for any $w_{0}\in U_{\delta }(u_{\lambda })$ the solution $%
v(t,w_{0})$ belongs to $U_{\varepsilon }(u_{\lambda })$ for all $t>0$.
Indeed, suppose the opposite, then since $v(t,w_{0})\in
\cC((0,T),\rH_{0}^{1}(\Omega ))$ there exists $t_{0}>0$ such that $||u_{\lambda
}-v(t_{0},w_{0})||_{1}=\varepsilon $. This implies that 
\begin{equation*}
d_{\varepsilon }\leq E_{\lambda }(v(t_{0},w_{0})).
\end{equation*}%
On the other hand, by Lemma \ref{lemst2} we have $E_{\lambda
}(v(t_{0},w_{0}))\leq E_{\lambda }(w_{0})$. Thus by \eqref{eqst3} one gets 
\begin{equation*}
d_{\varepsilon }\leq E_{\lambda }(v(t_{0},w_{0}))\leq E_{\lambda
}(w_{0})<d_{\varepsilon }-\sigma .
\end{equation*}%
This contradiction proves the claim. \fin
\bigskip

{\sc Proof of (III) Theorem \ref{Th2}}\quad Assume $\rN\geq 3$, $(\alpha ,\beta )\in \mathcal{E}_{s}(\rN)$ and $\Omega $ is
a strictly star-shaped domain with respect to the origin. By Corollary 15
from \cite{IlEg} it follows that there exists $\lambda ^{\ast }>0$ such that %
\eqref{1} has a flat ground state $u_{\lambda ^{\ast }}$ which $u_{\lambda
^{\ast }}\geq 0$ and $u_{\lambda ^{\ast }}\in \cC^{1,\gamma }(\overline{\Omega 
})\cap\cC^{2}(\Omega )$ for some $\gamma \in (0,1)$. Now applying (2),
Theorem \ref{Th1} we conclude that $u_{\lambda ^{\ast }}$ is a stable
non-negative stationary solution of the parabolic problem~\eqref{p1}.~\fin

\begin{rem}
\label{Remark 7.1} Related linearized stability results were obtained in 
\cite{BR} in working in Sobolev spaces in the framework of degenerate parabolic
equations of porous media type.
\end{rem}
\section{Linearized unstability}

In this Section, we prove statements \textbf{(I)} and \textbf{(II)} of
Theorem \ref{Th2}.

\begin{lemma}
\label{lemUn} Let $u_\lambda$ be a nonnegative weak solution of \eqref{1}
such that $E^{\prime \prime }(u_\lambda)<0$ then $u_{\lambda }$ is unstable
stationary solution of \eqref{p1} in the sense that $\lambda _{1}(-\Delta
-\lambda \beta u_\lambda^{\beta -1}+\alpha u_\lambda^{\alpha -1} )<0$.
\end{lemma}

\proof Let $u_{\lambda }$ be a nonnegative weak solution
of \ $SP(\alpha ,\beta ,\lambda )$. Then the corresponding linearized
problem at $u_{\lambda }$ is 
\begin{equation}
\left\{ 
\begin{array}{ll}
-\Delta \psi -(\lambda \beta u_{\lambda }^{\beta -1}-\alpha u_{\lambda
}^{\alpha -1})\psi =\mu \psi & \text{in }\Omega , \\ 
\,\psi =0 & \text{on }\partial \Omega .%
\end{array}%
\right.  \label{linearised}
\end{equation}%
Then there is a first eigenvalue $\mu _{1}$ to (\ref{linearised}) with a positive eigenfunction $%
\psi _{1}>0$ such that $\psi _{1}\in \cC^{2}(\Omega )\cap \cC_{0}^{1}(\overline{%
\Omega })$. The existence of $\mu _{1}$ is a particular case of the results
in \cite{DHMaagli} using the estimates on the boundary behavior of $%
u_{\lambda }$ obtained in \cite{D ambiguity}, \cite{D ambiguity 2}, namely
that

\begin{equation}
\underline{K}d(x)^{2/(1-\alpha )}\leq u_{\lambda }(x)\leq \overline{K}%
d(x)^{2/(1-\alpha )}\quad \hbox{for any
$x\in \overline{\Omega}$},  \label{near boundary}
\end{equation}%
for some constants $\overline{K}>\underline{K}>0.$ We shall sketch the
argument for the reader's convenience. From this estimates it follows that,
roughly speaking $u_{\lambda }(x)^{\alpha -1}$ "behaves like" $d(x)^{-2}$
and $u_{\lambda }(x)^{\beta -1}$ as $d(x)^{-2(1-\beta )/(1-\alpha )}$ with $%
\gamma :=2(1-\beta )/(1-\alpha )<2$ from $\alpha <\beta $. Then from the
used monotonicity properties of eigenvalues it is enough to show that a
first eigenvalue of the problem 
\begin{equation}
\left\{ 
\begin{array}{ll}
\displaystyle -\Delta w+\frac{\alpha }{d(x)^{2}}w-\frac{\lambda \beta }{d(x)^{\gamma }}%
w=\mu w & \text{in }\Omega , \\ 
\,w=0 & \text{on }\partial \Omega ,%
\end{array}%
\right.  \label{Hardy potential eigen}
\end{equation}%
is well-defined and has the usual properties. This is carried by reducing
the problem to an equivalent "fixed point" argument for an associated
(linear) eigenvalue problem. \ Assume first that $\mu >0$. Then (\ref{Hardy
potential eigen}) is equivalent to the existence of $\mu $ such that $r(\mu
)=1$, where $r(\mu )$ is the first eigenvalue for the associated problem 
\begin{equation}
\left\{ 
\begin{array}{lr}
\displaystyle -\Delta w+\frac{\alpha }{d(x)^{2}}w=r\left (\frac{\lambda \beta }{d(x)^{\gamma }}%
w+\mu w\right) & \text{in }\Omega , \\ 
\,w=0 & \text{on }\partial \Omega .%
\end{array}%
\right.  \label{Hardy erre}
\end{equation}%
That $r(\mu )>0$ is well-defined follows by showing that (\ref{Hardy erre})
is equivalently formulated as $Tw=rw$, with $T=i\circ P\circ F$, where $%
F:L^{2}(\Omega ,d^{\gamma })\rightarrow \rH^{-1}(\Omega )$ defined by 
\begin{equation*}
F(w)=\frac{\lambda \beta }{d(x)^{\gamma }}w+\mu w,
\end{equation*}

\noindent $P:\rH^{-1}(\Omega )\rightarrow \rH_{0}^{1}(\Omega )$ is the solution
operator for the linear problem 
\begin{equation}
\left\{ 
\begin{array}{lr}
\displaystyle -\Delta z+\frac{\alpha }{d(x)^{2}}z=h(x) & \text{in }\Omega , \\ 
\,w=0 & \text{on }\partial \Omega ,%
\end{array}%
\right.
\end{equation}%
for $h\in \rH^{-1}(\Omega )$, and $i:\rH_{0}^{1}(\Omega )\rightarrow
L^{2}(\Omega ,d^{\gamma })$ is the standard embedding. It is possible to
prove that $F$ and $P$ are continuous and $i$ is compact by using Hardy's
inequality and the Lax-Milgram Lemma (see\cite{BR},~\cite{DHMaagli}). Since 
$T$ is an irreductible compact linear operator and applying the weak maximum
principle, it is possible to apply Krein-Rutman's theorem in the formulation
in \cite{Daners KochMedina}. We have the variational formulation%
\begin{equation}
 r(\mu )=\inf_{w\in \rH_{0}^{1}(\Omega )\setminus \{0\}}\dfrac{\displaystyle\int_{\Omega
}\left (|\nabla w|^{2}+\frac{\alpha }{d(x)^{2}}w^{2}\right )\,dx}{\displaystyle  \lambda \beta
\int_{\Omega }\frac{w^{2}}{d(x)^{\gamma }}\,dx+\mu \int_{\Omega }w^{2}dx}.
\end{equation}%
Hence a positive eigenvalue exits if and only if there is a $\mu >0$ such
that $r(\mu )=1$. A completely analogous argument gives the formulation for $%
\mu <0$, namely with 
\begin{equation}
r_{1}(\mu )=\inf_{w\in \rH_{0}^{1}(\Omega )\setminus \{0\}}\frac{\displaystyle\int_{\Omega
}\left (|\nabla w|^{2}+\dfrac{\alpha }{d(x)^{2}}w^{2}+\mu w^{2}\right )\,dx}{\displaystyle\lambda \beta
\int_{\Omega }\frac{w^{2}}{d(x)^{\gamma }}\,dx}.
\end{equation}

\noindent Notice that $r(\mu )$ (resp. $r_{1}(\mu )$) is decreasing (resp.
increasing) in $\mu $. Then 
\begin{equation*}
r(0)=r_{1}(0)=\inf_{w\in \rH_{0}^{1}(\Omega )\setminus \{0\}}\frac{\displaystyle\int_{\Omega
}\left (|\nabla w|^{2}+\frac{\alpha }{d(x)^{2}}w^{2}\right )\,dx}{\displaystyle\lambda \beta
\int_{\Omega }\frac{w^{2}}{d(x)^{\gamma }}\,dx},
\end{equation*}%
and there exists a positive eigenvalue if $r(0)>1$ and a negative one if $%
r(0)<1.$

\noindent Coming back to our instability analysis, by Courant minimax
principle we have 
\begin{equation}
\mu _{1}=\inf_{\psi \in \rH_{0}^{1}(\Omega )\setminus \{0\}}\frac{\displaystyle \int_{\Omega
}\left (|\nabla \psi |^{2}-(\lambda \beta u_{\lambda }^{\beta -1}-\alpha
u_{\lambda }^{\alpha -1})\psi ^{2}\right )\,dx}{\displaystyle\int_{\Omega }|\psi |^{2}\,dx}
\label{Courant}
\end{equation}%
Let us put $\psi =u_{\lambda }$ in the minimizing functional of %
\eqref{Courant}. Then we get 
\begin{equation*}
\frac{\displaystyle\int_{\Omega }\left (|\nabla u_{\lambda }|^{2}-(\lambda \beta u_{\lambda
}^{\beta -1}-\alpha u_{\lambda }^{\alpha -1})u_{\lambda }^{2}\right )\,dx}{%
\displaystyle\int_{\Omega }|u_{\lambda }|^{2}\,dx}=\frac{E_{\lambda }^{\prime \prime
}(u_{\lambda })}{\displaystyle\int_{\Omega }|u_{\lambda }|^{2}\,dx}<0
\end{equation*}%
since by the assumption $E^{\prime \prime }(u_{\lambda })<0$. This yields by
the definition \eqref{Courant} that $\lambda _{1}(-\Delta -\lambda \beta
u_{\lambda }^{\beta -1}+\alpha u_{\lambda }^{\alpha -1}):=\mu _{1}<0$. Thus
we get an instability. \fin

\bigskip

{\sc Proof of (I), (II) Theorem \ref{Th2}}
\par
\noindent
{\sc Proof (I)}. \quad Assume $\rN=1,2$ and $(\alpha ,\beta )\in 
\mathcal{E}$. Let $u_{\lambda }$ be a free boundary solution of \eqref{1}.
Then since $\mathcal{E}=\mathcal{E}_{u}(\rN)$ statement 2) Lemma %
\ref{lemG} implies that $E_{\lambda }^{\prime \prime }(u_{\lambda })<0$.
However, this yields by Lemma \ref{lemUn} that $u_{\lambda }$ is a
linearized unstable stationary solution of the parabolic problem \eqref{p1}.
\par
\noindent 
{\sc Proof (II).}\quad Assume $\rN\geq 3$ and $(\alpha
,\beta )\in \mathcal{E}_{u}(\rN)$. Let $u_{\lambda }$ be a free boundary
solution of \eqref{1}. Then by 2), Lemma \ref{lemG} we have $E_{\lambda
}^{\prime \prime }(u_{\lambda })<0$. This yields as above by Lemma \ref%
{lemUn} that $u_{\lambda }$ is a linearized unstable stationary solution of
the parabolic problem \eqref{p1}.\fin


\section{Globally unstable ground state of \eqref{p1} in case $\protect\beta %
=1$}

In this Section, we prove statement \textbf{(2)}, Theorem \ref{Th2.2}.

Let us introduce the so called \textit{\ exterior potential well} (see \cite%
{PainSatt}) 
\begin{equation}
\mathcal{W}:=\{u\in \rH_{0}^{1}(\Omega ):E_{\lambda }(u)<\widehat{E}_{\lambda
},~~E_{\lambda }^{\prime }(u)<0\}.
\end{equation}

The proof of the theorem will be obtained from

\begin{lemma}
\label{lunst} If $v_{0}\in \mathcal{W}$, then $||v(t,v_{0})||_{L^{2}(\Omega
)}\rightarrow \infty $ as $t\rightarrow +\infty $.
\end{lemma}

\proof First we show that $\mathcal{W}$ is invariant
under the flow \eqref{p1}. Let $v(t,v_{0})$ be a weak solution of \eqref{p1}%
. Then using the additional regularity obtained in Section 2 we have 
\begin{equation*}
E_{\lambda }(v(t))\leq \int_{0}^{t}||v_{t}||_{L^{2}}^{2}ds+E_{\lambda
}(v(t))\leq E_{\lambda }(v_{0})<\widehat{E}_{\lambda }.
\end{equation*}%
for all $t>0$. Thus $v(t)$ may leave $\mathcal{W}$ only if there is a time $%
t_{0}>0$ such that $r_{\lambda }(v(t_{0}))=1$ (since, formally, $E_{\lambda
}^{\prime }(v(t_{0}))=0$). But then, by \eqref{maxr}, we have 
\begin{equation*}
E_{\lambda }(v(t_{0}))=\max_{r>0}E_{\lambda }(rv(t_{0}))\geq \widehat{E}_{\lambda
}.
\end{equation*}%
Thus we get a contradiction and indeed 
\begin{equation}
E_{\lambda }(v(t,v_{0}))<\widehat{E}_{\lambda },\quad E_{\lambda }^{\prime
}(v(t,v_{0}))<0~~~\forall t>0  \label{n1}
\end{equation}%
for any $v_{0}\in \mathcal{W}$. \fin

\bigskip

Furthermore, we have

\begin{prop}
\label{PrB} Assume that $v\in \rL^{\infty }(0,+\infty :\rH_{0}^{1}(\Omega ))$.
Then there exists $c_{0}<0,$ which does not depend on $t>0$ such that 
\begin{equation}
E_{\lambda }^{\prime }(v(t))\leq c_{0}<0~~~~~\text{for a.e. }t>0.
\label{nc0}
\end{equation}
\end{prop}

\proof By regularizing $v_{0}$ we can assume that $%
E_{\lambda }^{\prime }(v(t))$ is continuous in $t.$ Suppose, contrary to our
claim, that there is $(t_{m})$ such that the sequence $v_{m}:=v(t_{m})$, $%
m=1,2,...$ satisfies 
\begin{equation}
E_{\lambda }^{\prime }(v_{m})\rightarrow 0~~\mbox{as}~~m\rightarrow \infty .
\label{n2}
\end{equation}%
Note that by \eqref{n1} we have 
\begin{equation}
E_{\lambda }(v_{m})<\widehat{E}_{\lambda }~~\mbox{for}~~m=1,2,....  \label{n222}
\end{equation}%
By assumption $(v_{m})$ is bounded in $\rH_{0}^{1}(\Omega )$. Therefore we
have there are the following convergences (up choosing a subsequence) 
\begin{align}
& v_{m}\rightarrow \bar{v}~~\mbox{as}~~m\rightarrow \infty ~~\mbox{in}%
~~\rL^{p},~~1<p<2^{\ast }  \label{nc20} \\
& v_{m}\rightharpoondown \bar{v}~~\mbox{as}~~m\rightarrow \infty ~~%
\mbox{weakly in}~~\rH_{0}^{1}(\Omega )  \label{nc2} \\
& \lim_{m\rightarrow \infty }E_{\lambda }(v_{m})=a  \label{nn1}
\end{align}%
for some $\bar{v}\in \rH_{0}^{1}(\Omega )$ and $a\in \mathbb{R}$. Hence by the
weakly lower semi-continuity of $T(u)$ in $\rH_{0}^{1}(\Omega )$ we have 
\begin{align}
& E_{\lambda }(\bar{v})\leq \lim_{m\rightarrow \infty }E_{\lambda }(v_{m})=a
\label{n33} \\
& E_{\lambda }^{\prime }(\bar{v})\leq \lim_{m\rightarrow \infty }E_{\lambda
}^{\prime }(v_{m})=0.  \label{n3}
\end{align}%
Since $v\in \cC([0,T]:\rH_{0}^{1}(\Omega ))$ then by  Proposition \ref{proposition2.1} we have 
\begin{equation}
\int_{0}^{t}||v_{t}||_{L^{2}}^{2}ds+E_{\lambda }(v(t))\leq E_{\lambda
}(v(0)).  \label{ener}
\end{equation}%
Hence 
\begin{equation*}
a=\lim_{m\rightarrow \infty }E_{\lambda }(v_{m})\leq E_{\lambda }(v_{0})<%
\widehat{E}_{\lambda }
\end{equation*}%
for any $v_{0}\in \mathcal{W}$ and therefore $E_{\lambda }(\bar{v})<\widehat{E}%
_{\lambda }$. Observe that this implies a contradiction in case equality
holds in \eqref{n3}. Indeed, if $E_{\lambda }^{\prime }(\bar{v})=0$ then $r(%
\bar{v})=1$ and therefore \eqref{t}, \eqref{JL} and \eqref{minJomega} yield $%
E_{\lambda }(\bar{v})\geq \widehat{E}_{\lambda }$.

Suppose that $E_{\lambda }^{\prime }(\bar{v})<0$. Then there is $r\in (0,1)$
such that $E_{\lambda }^{\prime }(r\bar{v})=0$. Observe that \eqref{nc20}
and \eqref{nn1} imply 
\begin{equation}
\frac{1}{2}\lim_{m\rightarrow \infty }\rH_{\lambda }(v_{m})=a-\frac{1}{%
1+\alpha }A(\bar{v})
\end{equation}%
and \eqref{n2} implies 
\begin{equation}
\lim_{m\rightarrow \infty }\rH_{\lambda }(v_{m})=-A(\bar{v}).
\end{equation}%
From here we obtain 
$$
\begin{array}{ll}
E_{\lambda }(r\bar{v})& \displaystyle =\frac{r^{2}}{2}H_{\lambda }(\bar{v})+\frac{%
r^{1+\alpha }}{1+\alpha }A(\bar{v}) \\ [.3cm]
&  \displaystyle  \leq \frac{r^{2}}{2}\lim_{m\rightarrow \infty }\rH_{\lambda }(v_{m})+\frac{%
r^{1+\alpha }}{1+\alpha }A(\bar{v})\\ [.3cm]
& \displaystyle =\frac{1}{2}\lim_{m\rightarrow \infty
}\rH_{\lambda }(v_{m})+ \frac{1}{2}(r^{2}-1)\lim_{m\rightarrow \infty }\rH_{\lambda }(v_{m})+\frac{%
r^{1+\alpha }}{1+\alpha }A(\bar{v})\\ [.3cm]
&  \displaystyle =a-\frac{1}{1+\alpha }A(\bar{v})-\frac{1}{2}(r^{2}-1)A(\bar{v})+\frac{%
r^{1+\alpha }}{1+\alpha }A(\bar{v})\\ [.3cm]
&  \displaystyle =a+\left [-\frac{1}{1+\alpha }-\frac{1}{2}(r^{2}-1)+\frac{r^{1+\alpha }}{1+\alpha 
}\right ]A(\bar{v})
\end{array}%
$$
It is easy to see that 
\begin{equation*}
\max_{1\leq r\leq 1}\left \{\left [-\frac{1}{1+\alpha }-\frac{1}{2}(r^{2}-1)+\frac{%
r^{1+\alpha }}{1+\alpha }\right ]\right \}=0.
\end{equation*}%
Thus we get that $E_{\lambda }(r\bar{v})\leq a<\widehat{E}_{\lambda }$. However
this contradicts the definition of $\hat{E}_{\lambda },$ since $E_{\lambda
}^{\prime }(r\bar{v})=0$. This completes the proof of the proposition. \fin

\bigskip

Let us now conclude the proof of the Lemma. Suppose, contrary to our claim, that the set $(v(t))$, $t>0$ is bounded in $%
\rL^{2}(\Omega )$. Then this set is bounded also in $\rH_{0}^{1}(\Omega ),$
since $\rH_{\lambda }(v(t)):=T(v(t))-\lambda G(v(t))<0$ for all $t>0$.

Let us consider 
\begin{equation*}
y(t):=||v(t)||_{\rL^{2}}^{2},~~t\geq 0,
\end{equation*}%
where $v(t):=v(t,v_{0})$. Observe that 
\begin{equation*}
||v(t)||_{\rL^{2}}^{2}=||v_{0}||_{\rL^{2}}^{2}+2\int_{0}^{t}(v_{t}(s),v(s))\,ds
\end{equation*}%
and by \eqref{p1} 
\begin{equation*}
(v_{t}(s),v(s))=(\Delta v(s)+\lambda v(s)-|v(s)|^{\alpha
-1}v(s),v(s))=-E_{\lambda }^{\prime }(v(s)).
\end{equation*}%
Therefore 
\begin{equation}
y(t)=||v_{0}||_{\rL^{2}}^{2}-2\int_{0}^{t}E_{\lambda }^{\prime }(v(s))ds.
\end{equation}%
and 
\begin{equation*}
\frac{d}{dt}y(t)\equiv \dot{y}(t)=-2E_{\lambda }^{\prime }(v(t)).
\end{equation*}%
Hence estimates \eqref{nc0} of Proposition \ref{PrB} yields $\dot{y}%
(t)>-2c_{0}>0$ for all $t>0$ and therefore $y(t)=||v(t)||_{\rL^{2}}^{2}%
\rightarrow +\infty $ as $t\rightarrow \infty $. This completes the proof of
Lemma \ref{lunst}. \fin \newline
{\sc Conclusion of the proof of \textbf{(2)}, Theorem \ref{Th2.2}}\quad Let $u_{\lambda }$ be a ground state of \eqref{p1} and given any $\delta >0$%
. Observe that for any $r>1$ 
\begin{equation*}
E_{\lambda }(ru_{\lambda })<\widehat{E}_{\lambda }~~\mbox{and}~~E_{\lambda
}^{\prime }(ru_{\lambda })<0.
\end{equation*}%
Thus $ru_{\lambda }\in \mathcal{W}$ for any $r>1$ and by Lemma \ref{lunst} $%
||v(t;v_{0})||_{L^{2}}\rightarrow +\infty $ with $v_{0}=ru_{\lambda }$.
Therefore 
\begin{equation*}
||u_{\lambda }-v(t;v_{0})||_{\rL^{2}}\rightarrow +\infty ~~~\mbox{as}%
~~~t\rightarrow \infty .
\end{equation*}%
On the other hand, evidently $||u_{\lambda }-ru_{\lambda }||_{\rL^{2}}<\delta $
for sufficiently small $|r-1|$. This concludes the proof of Theorem \ref%
{Th2.2} \fin

\section*{Appendix. Existence of a ground state solution of \eqref{1R}}
\label{sec:appendix}
In this section we prove Lemma \ref{lem:r1}.

Consider 
\begin{equation}
\widehat{E}^{\infty }=\min \{J(v):~~v\in\rH_{0}^{1}(\Omega )\setminus \{0\},~\rH(v)<0\}.
\label{minJ}
\end{equation}

\begin{lemma}
\label{lemJ} There exists a minimizer $v$ of \eqref{minJ}.
\end{lemma}

\proof Let $(v_{m})$ be a minimizing sequence of %
\eqref{minJ}. Since $J(u)$ is a zero-homogeneous functional, we may assume
that $||v_{m}||_{1}=1$, $m=1,2,...,$. This implies that 
\begin{equation}
|\rH(v_{m})|<C<\infty ~~\mbox{uniformly on}~~m=1,2,....  \label{b1}
\end{equation}%
Observe that 
\begin{equation}
||v_{m}||_{\rL^{2}(\mathbb{R}^{\rN})}^{2}\equiv \int |v_{m}|^{2}dx>c_{1}>0
\label{b3}
\end{equation}%
uniformly on $m=1,2,...$. Indeed, if we suppose the contrary $\displaystyle \int
|v_{m}|^{2}dx\rightarrow 0$ as $m\rightarrow \infty $, then the assumption $%
||v_{m}||_{1}=1$, $m=1,2,...$ implies that $\displaystyle \int |\nabla
v_{m}|^{2}\,dx\rightarrow 1$ and therefore $\displaystyle \rH(v_{m})=\int |\nabla
v_{m}|^{2}\,dx-\int |v_{m}|^{2}\,dx\rightarrow 1$ as $m\rightarrow \infty $.
But this is impossible, since by the construction $\rH(v_{m})<0$.

Let us show that 
\begin{equation}
A(v_{m})>c_{0}>0~~\mbox{uniformly on}~~m=1,2,....  \label{b4}
\end{equation}%
Assume the opposite, that $A(v_{m})\rightarrow 0$ as $m\rightarrow \infty $.
Then $\displaystyle \int |v_{m}|^{2}dx\rightarrow 0$ as $m\rightarrow \infty $, since by H%
\"{o}lder and Sobolev inequalities 
\begin{equation*}
\int |v_{m}|^{2}dx\leq \left (\int |v_{m}|^{\alpha +1}dx\right )^{\frac{\kappa }{\alpha +1%
}}\left (\int |v_{m}|^{2^{\ast }}dx\right )^{\frac{\alpha +1-\kappa }{\alpha +1}}\leq
C_{0}A(v_{m})^{\frac{\kappa }{\alpha +1}}||v_{m}||_{1}^{2^{\ast }\frac{%
\alpha +1-\kappa }{\alpha +1}},
\end{equation*}%
where $\kappa =\dfrac{(\alpha +1)(2^{\ast }-2)}{2^{\ast }-\alpha +1}$. But
this contradicts \eqref{b3}.

Observe that \eqref{fibF}, \eqref{b1} and \eqref{b4} yield 
\begin{equation}  \label{b0}
\widehat{E}^\infty >0,
\end{equation}
and we have 
\begin{equation}  \label{b5}
0<c_0<||v_m||_{\rL^{1+\alpha}}^{1+\alpha} \equiv A(v_m) <C_1<+\infty
\end{equation}
uniformly on $m=1,2,\ldots $

We need the following lemma \cite[Lemma I.1, p.231]{lions}

\begin{lemma}
\label{Lions} Let $1\leq q<+\infty $ with $q\leq 2^{\ast }$ if $\rN\geq 3$.
Assume that $(w_{n})$ is bounded in $\rH_{0}^{1}(\mathbb{R}^{\rN})$ and $\rL^{q}(%
\mathbb{R}^{N})$, and 
\begin{equation*}
\sup_{y\in \mathbb{R}^{N}}\int_{y+B_{R}}|w_{n}|^{q}dx\rightarrow 0~~\mbox{as}%
~~n\rightarrow \infty ,~~\mbox{for some}~~R>0
\end{equation*}%
Then $||w_{n}||_{\rL^{\beta }}\rightarrow 0$ for $\beta \in (q,2^{\ast })$.
\end{lemma}

Let $R>0$. Observe that 
\begin{equation}
\liminf_{m\rightarrow \infty }\sup_{y\in \mathbb{R}^{\rN}}%
\int_{y+B_{R}}|v_{m}|^{1+\alpha }dx:=\delta >0.
\end{equation}%
Indeed, let us assume that 
\begin{equation*}
\liminf_{m\rightarrow \infty }\sup_{y\in \mathbb{R}^{\rN}}%
\int_{y+B_{R}}|v_{m}|^{1+\alpha }dx=0.
\end{equation*}%
Then by Lemma \ref{Lions} we have $||v_{m}||_{\rL^{2}}\rightarrow 0$ as $%
m\rightarrow \infty $. But this contradicts \eqref{b3}.

Thus there is a sequence $\{y_{m}\}\subset \mathbb{R}^{\rN}$ such that 
\begin{equation*}
\int_{y_{m}+B_{R}}|v_{m}|^{1+\alpha }dx>\frac{\delta }{2},~~~m=1,2,....
\end{equation*}%
Introduce $u_{m}:=v_{m}(\cdot +y_{m})$, $m=1,2,...$. Then 
\begin{equation}
\int_{B_{R}}|u_{m}|^{1+\alpha }dx>\frac{\delta }{2},~~~m=1,2,...,  \label{nz}
\end{equation}%
and $\{u_{m}\}$ is a minimizing sequence of \eqref{minJ}.

Furthermore, by the zero-homogeneity of $J(u)$ now we may normalize the
sequence $\{u_{m}\}$ (again denoted by $\{u_{m}\}$) such that 
\begin{equation}
A(u_{m})=1,~~~m=1,2,....
\end{equation}%
Then \eqref{b5} implies that the renormalized sequence $\{u_{m}\}$ will be
again bounded in $\rH^{1}(\mathbb{R}^{\rN})$. Thus by Eberlein-Smulian theorem
there is a subsequence of $\{u_{m}\}$ (again denoting $\{u_{m}\}$) and a limit
point $\bar{u}\in \rH_{0}^{1}(\Omega )$ such that 
\begin{equation}
u_{m}\rightharpoondown \bar{u}~~\mbox{weakly in }~~\rH_{0}^{1}(\Omega )~~%
\mbox{as}~~m\rightarrow \infty .  \label{con1}
\end{equation}%
Furthermore 
\begin{equation}
u_{m}\rightarrow \bar{u}~~\mbox{a.e. on }~~\mathbb{R}^{\rN}~~\mbox{as}%
~~m\rightarrow \infty ,  \label{6}
\end{equation}%
and for $2<q<2^{\ast }$ 
\begin{equation}
u_{m}\rightarrow \bar{u}~~\mbox{in}~~\rL_{loc}^{q}~~\mbox{as}~~m\rightarrow
\infty ,
\end{equation}%
since by Rellich-Kondrachov theorem $\rH_{0}^{1}(B_{R})$ is compactly embedded
in $\rL^{q}(B_{R})$ for $2<q<2^{\ast }$ and any $B_{R}:=\{x\in \mathbb{R}%
^{\rN}:|x|\leq R\}$, $R>0$. Note that \eqref{nz} implies that 
\begin{equation*}
\bar{u}\neq 0.
\end{equation*}%
We need the Brezis-Lieb lemma \cite{BrezisLieb}:

\begin{lemma}
Let $\Omega $ be an open subset of $\mathbb{R}^{\rN}$ and let $\{w_{n}\}\subset
\rL^{q}(\Omega )$, $1\leq q<\infty $. If

\begin{description}
\item[a)] $\{w_{n}\}$ bounded in $\rL^{q}(\Omega )$,

\item[b)] $w_{n}\rightarrow w$ a.e. on $\Omega $, then 
\begin{equation*}
\lim_{n\rightarrow \infty
}\big (||w_{n}||_{\rL^{q}}^{q}-||w_{n}-w||_{\rL^{q}}^{q}\big)=||w||_{\rL^{q}}^{q}.
\end{equation*}
\end{description}
\end{lemma}

Let us denote $\omega _{m}:=u_{m}-\bar{u}$. Then Brezis-Lieb lemma yields 
\begin{equation}
1=A(\bar{u})+\lim_{m\rightarrow \infty }A(\omega _{m}).  \label{con2}
\end{equation}%
Observe 
\begin{equation}
\rH(\omega _{m})=\rH(\bar{u})+\rH(u_{m})-\rH^{\prime }(u_{m})(\bar{u}).  \label{b6}
\end{equation}%
Note that due to weak convergence \eqref{con1} we have $H^{\prime }(\omega
_{m})(u)\rightarrow 0$ as $m\rightarrow \infty $. Therefore, $H(\omega
_{m})<0$ for sufficiently large $m$, since $H(u)<0$ and $H(u_{m})<0$ for $%
m=1,2,...$. On the other hand 
\begin{equation*}
\rH(u_{m})=\rH(\bar{u})+\rH(\omega _{m})+\rH^{\prime }(\omega _{m})(\bar{u}),
\end{equation*}%
and therefore 
\begin{equation}
\lim_{m\rightarrow \infty }\rH(u_{m})=\rH(\bar{u})+\lim_{m\rightarrow \infty
}\rH(\omega _{m}).  \label{b7}
\end{equation}%
Observe that \eqref{minJ} implies that for any $v\in \rH_{0}^{1}(\Omega
)\setminus \{0\}$ s.t. $\rH(v)<0$ it holds 
\begin{equation}
-\rH(v)\leq k_{\alpha }\frac{A(v)^{\frac{2}{1+\alpha }}}{\widehat{E}^{\infty }}
\end{equation}%
where 
\begin{equation*}
k_{\alpha }=\left (\frac{(1-\alpha )}{2(1+\alpha )}\right )^{\frac{1-\alpha }{1+\alpha }}.
\end{equation*}%
Hence 
\begin{equation*}
-\rH(\bar{u})\leq k_{\alpha }\frac{A(\bar{u})^{\frac{2}{1+\alpha }}}{\widehat{E}%
^{\infty }}
\end{equation*}%
and 
\begin{equation}
-\rH(\omega _{m})\leq k_{\alpha }\frac{A(\omega _{m})^{\frac{2}{1+\alpha }}}{%
\widehat{E}^{\infty }},  \label{con3}
\end{equation}%
for sufficient large $m$. Since $A(u_{m})=1$, we have 
\begin{equation*}
\lim_{m\rightarrow \infty }k_{\alpha }\frac{1}{(-\rH(u_{m}))}=\widehat{E}^{\infty
}.
\end{equation*}%
Hence we have 
$$
\begin{array}{ll}
\displaystyle k_{\alpha }\frac{1}{\widehat{E}^{\infty }}& \displaystyle= \lim_{m\rightarrow \infty
}(-\rH(u_{m}))\\ [.3cm]
&\displaystyle =-\rH(\bar{u})+\lim_{m\rightarrow \infty }(-\rH(\omega _{m}))\\ [.3cm]
& \displaystyle\le k_{\alpha }\frac{A(\bar{u})^{\frac{2}{1+\alpha }}}{\widehat{E}^{\infty }}%
+\lim_{m\rightarrow \infty }k_{\alpha }\frac{A(\omega _{m})^{\frac{2}{%
1+\alpha }}}{\widehat{E}^{\infty }}\\ [.3cm]
&\displaystyle=k_{\alpha }\frac{1}{\widehat{E}^{\infty }}\left(
A(\bar{u})^{\frac{2}{1+\alpha }}+(1-A(\bar{u}))^{\frac{2}{1+\alpha }}\right).
\end{array}%
$$
Note since $\frac{2}{1+\alpha }>1$, then $f(r):=r^{\frac{2}{1+\alpha }%
}+(1-r)^{\frac{2}{1+\alpha }}\geq 1$ for $r\in \lbrack 0,1]$ and $f(r)=1$
iff $r=0$ or $r=1$. Thus we have 
\begin{equation*}
A(\bar{u})=1~~\mbox{or}~~A(\bar{u})=0.
\end{equation*}%
Now taking into account that $\bar{u}\neq 0$ we get that $A(\bar{u})=1$.
Hence by \eqref{con2} we obtain $A(\omega _{m})\rightarrow 0$ as $%
m\rightarrow \infty $ and consequently by \eqref{con3} we have $(-\rH(\omega
_{m}))\rightarrow 0$ as $m\rightarrow \infty $. From here it is not hard to
conclude that $u_{m}\rightarrow \bar{u}$ strongly in $\rH^{1}(\mathbb{R}^{\rN})$
and therefore $J(\bar{u})=\widehat{E}^{\infty }$. Thus $\bar{u}$ is a minimizer
of \eqref{minJ}.~\fin
\medskip
\par
\noindent 
{\sc Proof of Lemma \ref{lem:r1}}.\quad By Lemma \ref{lemJ} there
exists a minimizer $\bar{u}$ of \eqref{minJ}. Since $J$ is an even
functional then $|\bar{u}|$ is also a minimizer of \eqref{minJ}. Thus we may
assume that $\overline{u}$ is nonnegative function. By Proposition \ref%
{prop:fc} it follows that $u=r(\bar{u})\bar{u}$ is a weak solution of %
\eqref{1R} which is nonnegative since $r(\bar{u})>0$. By regularity theory
we derive that $u\in \cC^{2}(\mathbb{R}^{\rN})$.\fin

\flushright{
\begin{tabular}{lll}
J.I. D\'{\i}az & J. Hern\'{a}ndez  & Y. Il'yasov \\
Instituto de Matem\'{a}tica Interdisciplinar & Departamento de Matem\'{a}ticas & Institute of Mathematics\\
Universidad Complutense de Madrid &  Universidad Aut\'{o}noma de Madrid &  Ufa Science Center of RAS \\
28040 Madrid, Spain & 28049 Cantoblanco, Madrid, Spain & Chernyshevsky Str.,\\& & Ufa, Russia \\
{\tt jidiaz@ucm.es} & {\tt jesus.hernandez@uam.es} & {\tt ilyasov02@gmail.com}
\end{tabular}
}

\end{document}